\documentclass[a4paper,12pt,oneside]{article} 

\title{Boundedness of Log Canonical Surface Generalized Polarized Pairs}
\author{Stefano Filipazzi}
\date{\vspace{-5ex}} 

\usepackage{amsthm}
\usepackage{amssymb}
\usepackage{amsmath,amscd}
\usepackage[english,italian]{babel}
\usepackage{latexsym}
\usepackage{graphicx,subfigure}
\usepackage[a4paper,top=30mm,bottom=30mm,left=30mm,right=30mm]{geometry}
\usepackage{tikz}
\usetikzlibrary{matrix,arrows,decorations.pathmorphing}

\usepackage{hyperref}

\usepackage{mathtools}

\usepackage[babel]{csquotes}
\usepackage[maxnames=99,style=alphabetic,backend=bibtex]{biblatex}
\usepackage{scalerel,stackengine}
\stackMath
\newcommand\reallywidehat[1]{%
\savestack{\tmpbox}{\stretchto{%
  \scaleto{%
    \scalerel*[\widthof{\ensuremath{#1}}]{\kern-.6pt\bigwedge\kern-.6pt}%
    {\rule[-\textheight/2]{1ex}{\textheight}}
  }{\textheight}%
}{0.5ex}}%
\stackon[1pt]{#1}{\tmpbox}%
}

\theoremstyle{definition}
\newtheorem{teo}{Theorem}[section]
\newtheorem{cor}[teo]{Corollary}
\newtheorem{lemma}[teo]{Lemma}
\newtheorem{prop}[teo]{Proposition}
\newtheorem{oss}[teo]{Observation}
\newtheorem{rem}[teo]{Remark}

\newtheorem*{teo*}{Theorem}
\newtheorem*{prop*}{Proposition}

\newtheorem*{def*}{Definition}
\newtheorem*{oss*}{Observation}
\newtheorem*{rem*}{Remark}

\newtheorem*{ex*}{Example}
\newtheorem*{cor*}{Corollary}

\newtheorem*{ack*}{Acknowledgements}

\newtheorem*{dimo1*}{Proof}
\newtheorem{thmx}{Theorem}
\newtheorem{propx}[thmx]{Proposition}

\newcommand{\Ox}{\mathcal{O}}

\newcommand{\Q}{\mathbb{Q}}
\newcommand{\R}{\mathbb{R}}
\newcommand{\C}{\mathbb{C}}
\newcommand{\Z}{\mathbb{Z}}
\newcommand{\N}{\mathbb{N}}
\newcommand{\Pro}{\mathbb{P}}

\newcommand{\du}{\,\check{\vrule height1.3ex width0pt}}

\newcommand{\mult}{\mathrm{mult}}
\newcommand{\vol}{\mathrm{vol}}
\newcommand{\coeff}{\mathrm{coeff}}
\newcommand{\red}{\mathrm{red}}

\newenvironment{claim}[1]{\par\noindent\underline{Claim:}\space#1}{}

\usepackage{lipsum}

\newcommand{\Addresses}{{
  \bigskip
  \footnotesize

  S.~Filipazzi, \textsc{Department of Mathematics, University of Utah,
    Salt Lake City,\\ UT 84112, USA}\par\nopagebreak
  \textit{E-mail address}: \texttt{filipazz@math.utah.edu}

}}

\begin{document}

\selectlanguage{english}

\maketitle

\begin{abstract}

In this paper, we study the behavior of the sets of volumes of the form $\vol(X,K_X+B+M)$, where $(X,B)$ is a log canonical pair, and $M$ is a nef $\R$-divisor. After a first analysis of some general properties, we focus on the case when $M$ is $\Q$-Cartier with given Cartier index, and $B$ has coefficients in a given DCC set. First, we show that such sets of volumes satisfy the DCC property in the case of surfaces. Once this is established, we show that surface pairs with given volume and for which $K_X+B+M$ is ample form a log bounded family. These generalize results due to Alexeev \cite{Ale94}.

\end{abstract}

\tableofcontents

\section{Introduction}

The study of the volume of log pairs plays a crucial role in the understanding of varieties of log general type and the proof of their boundedness. In particular, we have the following result by Hacon, McKernan and Xu.

\begin{thmx}[{\cite[Theorem 1.9]{HMX13}}] \label{teoA}
Fix a set $\Lambda \subset [0,1]$ which satisfies the DCC. Let $\mathfrak{D}$ be a set of simple normal crossings pairs $(X,B)$, which is log birationally bounded, such that, if $(X,B) \in \mathfrak{D}$, then the coefficients of $B$ belong to $\Lambda$.

Then the set
\begin{equation*}
\lbrace \vol(X,K_X+B)|(X,B) \in \mathfrak{D} \rbrace
\end{equation*}
satisfies the DCC.
\end{thmx}

The above theorem is a crucial result in the proof of the boundedness of varieties of log general type. Recall that $\mathfrak{F}_{\mathrm{slc}}(n,\Lambda,d)$ denotes the set of projective semi-log canonical pairs $(X,B)$ of dimension $n$ such that $K_X+B$ is ample, $(K_X+B)^n=d$, and the coefficients of $B$ lie in $\Lambda \subset [0,1]$. Hacon, Mckernan and Xu proved the following.

\begin{thmx}[{\cite[Theorem 1.2.1]{HMX}}] \label{teoB}
Fix $n \in \N$, a set $\Lambda \subset [0,1] \cap \Q$ which satisfies the DCC and $d > 0$. Then the set $\mathfrak{F}_{\mathrm{slc}}(n,\Lambda,d)$ is bounded; that is, there is a projective morphism of quasi-projective varieties $\pi : \mathcal{X} \rightarrow T$ and a $\Q$-divisor $\mathcal{B}$ on $\mathcal{X}$ such that the set of pairs $\lbrace (\mathcal{X}_t,\mathcal{B}_t) | t \in T \rbrace$ given by the fibers of $\pi$ is in bijection with the elements of $\mathfrak{F}_{\mathrm{slc}}(n,\Lambda,d)$.
\end{thmx}

It is relevant to point out that versions for surfaces of Theorem \ref{teoA} and Theorem \ref{teoB} were first proved by Alexeev \cite{Ale94}.

Recently, Birkar and Zhang introduced the notion of generalized polarized pair \cite{BZ}. This kind of pair arises naturally in certain situations, such as the canonical bundle formula \cite{FM00}. Furthermore, these generalized polarized pairs play an important role in recent developments, such as the study of the Iitaka fibration \cite{BZ}, and the proof of the BAB conjecture \cite{Bir16a}, \cite{Bir16b}. Therefore, it is interesting to further investigate the properties of generalized polarized pairs.

In the hope that a program in the spirit of \cite{HMX} could be accomplished in the setting of generalized polarized pairs, we start investigating the properties of volumes in this new context. In particular, we are motivated by the following fact shown in \cite{BZ}.

\begin{propx}[{\cite{BZ}}] \label{propC}
Let $\Lambda \subset [0,1]$ be a DCC set, and let $d,r$ be positive integers. Then there exists a number $\beta>0$ depending only on $\Lambda,d,r$ such that, if
\begin{itemize}
\item $(X,B)$ is projective log canonical of dimension $d$;
\item the coefficients of $B$ are in $\Lambda$;
\item $rM$ is a nef Cartier divisor;
\item $K_X+B+M$ is a big divisor;
\end{itemize}
then $\vol(X,K_X+B+M) > \beta$.
\end{propx}

The above statement in not explicitly formulated in \cite{BZ}, and just a weaker version is presented in \cite[Step 7, Proposition 3.4]{BZ}. On the other hand, Proposition \ref{propC} is equivalent to the following.

\begin{thmx}[{\cite[Theorem 1.3]{BZ}}] \label{BirkarZhangBirational}
Let $\Lambda \subset [0,1]$ be a DCC set, and $d,r$ positive integers. Then there is a positive integer $m_0$ depending only on $\Lambda, d, r$ such that, if
\begin{itemize}
\item $(X,B)$ is a projective log canonical pair of dimension $d$;
\item the coefficients of $B$ are in $\Lambda$;
\item $rM$ is a nef Cartier divisor;
\item $K_X+B+M$ is big;
\end{itemize}
then the linear system $|m(K_X+B+M)|$ defines a birational map if $m$ is a positive integer divisible by $m_0$.
\end{thmx}

\begin{rem}
Given an $\R$-divisor $D$, by the linear series $|D|$ we denote the linear series of the divisor $\lfloor D \rfloor$, i.e. $|\lfloor D \rfloor |$.
\end{rem}

This hints that the above set of volumes of generalized polarized pairs could satisfy the DCC property. In this paper, we prove the following, which is a generalization of \cite[Theorem 8.2]{Ale94}.

\begin{teo} \label{MainResult}
Let $\Lambda \subset [0,1]$ be a DCC set, and let $r$ be a positive integer. Then the set
\begin{equation}
\mathfrak{V} \coloneqq \lbrace \vol(X,K_X+B+M)|(X,B)\; \mathrm{is \; lc \; surface,\;} B \in \Lambda,\; rM \mathrm{\; is \; nef \; and \; Cartier}   \rbrace
\end{equation}
satisfies the DCC property.
\end{teo}

The statement of Theorem \ref{MainResult} and the constructions performed in its proof lead us to a more interesting statement concerning boundedness of surfaces. Before presenting the related statement, which generalizes \cite[Theorem 9.2]{Ale94}, we need to introduce some terminology.

\begin{def*}
A \emph{generalized (semi) log canonical model} $(X',B'+M')$ is a projective generalized (semi) log canonical generalized polarized pair $(X',B'+M')$ such that $K_{X'}+B'+M'$ is ample. Fix positive integers $d$ and $r$, a positive number $w$ and a set $\Lambda \subset [0,1]$. Then, we denote by $\mathfrak{F}_{gslc}(d,w,\Lambda,r)$ the set of all $d$-dimensional generalized semi-log canonical models $(X',B'+M')$ such that the coefficients of $B'$ belong to $\Lambda$, $rM'$ is Cartier, $K_{X'}+B'+M'$ is $\Q$-Cartier and $(K_{X'}+B'+M')^d=w$.
\end{def*}

\begin{teo} \label{BigResult}
Fix $\Lambda \subset [0,1]\cap \Q$ a DCC set, a positive integer $r$ and let $w$ be a positive number. Then, the set $\mathfrak{F}_{gslc}(2,w,\Lambda,r)$ is bounded. That is, there is a positive integer $N=N(2,w,\Lambda,r)$ such that if $(X',B'+M') \in \mathfrak{F}_{gslc}(2,w,\Lambda,r)$, then $N(K_{X'}+B'+M')$ is very ample. In particular, the coefficients of $B'$ belong to a fixed finite set $\Lambda_0 \subset \Lambda$.
\end{teo}

\begin{rem}
In view of Proposition \ref{boundHB}, we conclude that the pairs $(X',B')$ appearing in Theorem \ref{BigResult} are bounded. Furthermore, since we can write $N(K_{X'} + B' + M') \sim H$ for some very ample divisor that deforms along the bounding family, we conclude that $M'$ is bounded as well up to $\Q$-linear equivalence.
\end{rem}

As a first step, we consider a special case of Theorem \ref{MainResult}: We study the behavior of the volume of higher models of a fixed generalized polarized pair. The techniques involved are a generalization of the ones used in \cite{HMX}.

\begin{teo} \label{DCCcase1}
Let $\Lambda \subset [0,1]$ be a DCC set, and let $(Z,D)$ be a simple normal crossing pair where $D$ is reduced, and $M$ is a fixed $\mathbb{R}$-Cartier $\mathbb{R}$-divisor on $Z$. Let $\mathfrak{D}$ be the set of all projective simple normal crossing pairs $(X,B)$ such that $\coeff (B) \subset \Lambda$, $f: X \rightarrow Z$ is a birational morphism and $f_*B \leq D$. Then, the set
\begin{equation}
\lbrace \vol(X,K_X+B+f^*M)|(X,B)\in \mathfrak{D} \rbrace
\end{equation}
satisfies the DCC.
\end{teo}

\begin{rem}
Since we are interested in volumes of the form $\vol(X,K_X+B+M)$ where $M$ is nef and $(X,B)$ is a log canonical pair, the simple normal crossing assumptions in the statements are not restrictive. Indeed, we are free to take a log resolution of $(X,B)$, take as new boundary the strict transform of $B$ plus the reduced exceptional divisor, and pull back $M$. This process preserves the volume, as shown in Lemma \ref{smoothModel}.
\end{rem}

Once a result concerning a single birational class is established, one would like to control volumes in families. By perturbing the nef part by a small ample divisor and using continuity of the volume function, we can apply the well known result for usual pairs \cite[Theorem 2.6.2]{HMX} and get a statement about deformation invariance of volumes.

\begin{teo} \label{DefInvariance}
Let $\mathcal{(X,\mathcal{B})}\rightarrow T$ be a log-smooth family, where $T\supset{\lbrace x_i \rbrace_{i \geq 1}}$. Denote by $(X_i,B_i)$ the log pair cut by $\mathcal{B}$ over $x_i$. Assume that
\begin{itemize}
\item $0 \leq \mathcal{B} \leq \mathcal{B}_\red$;
\item there is an $\R$-divisor $\mathcal{M}$ on $\mathcal{X}$ such that $M_i \coloneqq \mathcal{M}|_{X_i}$ is nef for every $i$;
\end{itemize}
Then, we have $\vol (X_i,K_{X_i}+B_i+M_i)=\vol (X_j,K_{X_j}+B_j+M_j)$ for every $i,j \in \N$.
\end{teo}

Once Theorem \ref{DCCcase1} and Theorem \ref{DefInvariance} are established, the proof strategy for Theorem \ref{MainResult} becomes clear. Arguing by contradiction, there is a sequence $\lbrace (X_i,B_i+M_i) \rbrace_{i \geq 1}$ as in the statement of Theorem \ref{MainResult} such that $\lbrace \vol(X_i,K_{X_i}+B_i+M_i) \rbrace_{i \geq 1}$ forms a strictly decreasing sequence. Using the techniques in \cite{HMX}, it follows that the pairs $\lbrace (X_i,B_i)\rbrace_{i \geq 1}$ are log birationally bounded. Let $(\mathcal{X},\mathcal{B})\rightarrow T$ be a log smooth bounding family, and $\mathcal{X}_i$ the birational model of $X_i$. Without loss of generality, the birational maps $g_i : X_i \dashrightarrow \mathcal{X}_i$ are actual morphisms. Assume that we can also produce a divisor $\mathcal{M}$ on $\mathcal{X}$ that bounds the numerical class of the $g_{i,*}M_i$'s, and $M_i = g^*_i g_{i,*}M_i$. Then, by Theorem \ref{DCCcase1}, the volumes of the birational models dominating the $\mathcal{X}_i$'s satisfy the DCC. Thus, deformation invariance of volumes allows us to compare the DCC sets living over different fibers and conclude that they agree, providing the required contradiction.

The key parts of the above argument are the construction of $\mathcal{M}$ and assuring that the condition $M_i = g^*_i g_{i,*}M_i$ holds. In the surface case, there are some further techniques that can be applied to achieve these goals. Relying on the fact that the intersection matrix of the exceptional curves in a birational morphism between normal surfaces is negative definite, we can produce a natural candidate for $\mathcal{M}$. Then, the second condition is related to the following questions. Given a morphism $f : X \rightarrow Y$ and a nef divisor $M$ on $X$, when does $M$ descend to $Y$? How high of an intermediate model carries all the information relative to $M$? In the case of surfaces, we prove the following effective statement.

\begin{teo} \label{structureSurfaceNef}
Let $f: Y \rightarrow X$ be a birational morphism of smooth surfaces. Let $M_Y$ and $M$ be nef Cartier divisors on $Y$ and $X$ respectively. Assume $f_*M_Y=M$, and that $M_Y^2=M^2-k$. Then, there exist a smooth surface $X'$, a nef divisor $M'$ on $X'$, morphisms $g:Y \rightarrow X'$ and $h:X' \rightarrow X$ such that:
\begin{itemize}
\item $f=h \circ g$;
\item $h_*M'=M$;
\item $g^*M'=M_Y$;
\item $X'$ is obtained from $X$ by blowing up at most $k$ points.
\end{itemize}
Furthermore, such an $X'$ is obtained by blowing-down all $f$-exceptional $(-1)$-curves until $M \cdot E > 0$ for all exceptional $(-1)$-curves.
\end{teo}

Theorem \ref{structureSurfaceNef} plays an essential role in the proof of Theorem \ref{MainResult}. In particular, it allows us to reduce to the case where we can apply Theorem \ref{DCCcase1}.

In order to prove Theorem \ref{BigResult}, we can reduce to the case when the semi-log canonical models are disjoint unions of log canonical ones. 
Indeed, once this partial result is established, the full statement of Theorem \ref{BigResult} follows by an effective basepoint-free theorem for semi-log canonical surfaces due to Fujino \cite{Fuj17}. Then, the key result is the following one.

\begin{teo} \label{BoundednessLCModels}
Fix $\Lambda \subset [0,1]\cap \Q$ a DCC set, a positive integer $r$ and let $w$ be a positive number. Let $\mathfrak{F}_{glc}(2,w,\Lambda,r)$ be the set of generalized polarized surface pairs $(X,B+M)$ that are disjoint union of generalized log canonical models $(X_j,B_j+M_j)$, where $rM$ is nef and Cartier, $\coeff(B_j) \subset \Lambda$ and $(K_X+B+M)^2=w$. Then $\mathfrak{F}_{glc}(2,w,\Lambda,r)$ is bounded. 
\end{teo}

\begin{rem} \label{reduction to surface pairs}
In general, given a generalized polarized pair $(X',B'+M')$, the divisor $K_{X'}+B'$ is not necessairily $\R$-Cartier. On the other hand, in the case of surfaces we can compute its discrepancies numerically \cite[p. 112]{KM}. If $(X',B'+M')$ is generalized log canonical, the Negativity Lemma \cite[Lemma 3.39]{KM} applied to $M$ and the numerical pullback of $M'$ implies that $(X',B')$ is numerically log canonical. Thus, $(X',B')$ is honestly log canonical, and then $M'$ is $\R$-Cartier. Then, Lemma \ref{push nef} and Remark \ref{push nef singular} imply that $M'$ is nef. Therefore, when dealing with surfaces, it is not restrictive to consider a log canonical pair $(X',B')$ decorated with a nef divisor $M$ on a birational model $X$ in place of a generalized polarized pair.

Although the statement of Theorem \ref{BoundednessLCModels} talks about boundedness for the base model of a generalized polarized pair, i.e. what we usually denote by $(X',B'+M')$, it is relevant to point out how this, together with Theorem \ref{structureSurfaceNef}, also implies boundedness for the nef part living on a higher birational model.
\end{rem}

The proof of Theorem \ref{BoundednessLCModels} builds on the constructions made for Theorem \ref{MainResult}. Again, we argue by contradiction. Thus, there is a sequence $\lbrace (X_i,B_i+M_i) \rbrace_{i \geq 1}$ as in Theorem \ref{BoundednessLCModels} for which there is no uniform $N$ such that $N(K_{X_i}+B_i+M_i)$ is very ample for all $i \geq 1$. By Theorem \ref{MainResult}, we may assume that all the $X_i$'s are irreducible. Then, by the proof of Theorem \ref{MainResult}, the sequence is log birationally bounded by a log smooth family $(\mathcal{X},\mathcal{B}+ \mathcal{M}) \rightarrow T$. By techniques developed in \cite{HMX}, we can reduce to the case where the birational maps $q_i : \mathcal{X}_i \dashrightarrow X_i$ are actual morphisms. Furthermore, we may assume that the $q_i$-exceptional divisors deform along certain components of $\mathcal{B}$ that are independent of $i$. This implies that the Cartier index of $X_i$ is uniform with respect to $i$. 
In particular, this guarantees the existence of a positive integer $l$ such that $l(K_{X_i}+B_i+M_i)$ is integral and Cartier for all $i$. Since $K_{X_i}+B_i+M_i$ is ample and $M_i$ is nef, we can then apply an effective basepoint-free theorem for semi-log canonical surfaces by Fujino \cite{Fuj17}. In particular, we show that $12l(K_{X_i}+B_i+M_i)$ is very ample, and  the required contradiction follows.

In the last section of the paper, we consider some examples of generalized polarized pairs that show how good properties of pairs do not, in general, extend to this new setup. First, we show how Theorem \ref{DefInvariance} is somehow optimal: Indeed, we provide an example of a family as in Theorem \ref{DefInvariance} for which deformation invariance of plurigenera does not hold. Then, we show that for a general fiber of such a family the generalized pluricanonical ring $R(Y,K_Y+B+M)$ is not finitely generated, although $(Y,B)$ is log canonical and $K_Y+B+M$ is nef and big.

\begin{ack*}
The author would like to thank his advisor Christopher D. Hacon for suggesting the problem, for his insightful suggestions and encouragement. He would also like to thank Karl Schwede and Roberto Svaldi for helpful conversations.
\end{ack*}

\section{General Properties}

In this section we introduce a few properties of generalized polarized pairs that hold in any dimension. Also, we recall some facts and fix our notation.

\subsection{Background and Notation}

In this paper we work over the field of complex numbers $\C$. In particular, all the constructions and statements have to be understood in this a setting. Now, we recall a few properties and definitions that will be relevant in the following.

First, we include the usual definition of boundary and log pair.

\begin{def*}
Let $X$ be a normal projective variety. A {\it boundary} $B$ is an effective $\R$-divisor such that $K_X+B$ is $\R$-Cartier. A {\it log pair} (or simply a {\it pair}) $(X,B)$ is the datum of a normal projective variety and a pair.
\end{def*}

Now that we have fixed our notation, we can introduce the definition of generalized polarized pair, and compare it with the usual notion of pair.

\begin{def*}
A {\it generalized polarized pair} consists of a normal variety $X'$, equipped with morphisms $X \xrightarrow{f} X' \rightarrow Z$, where $f$ is birational and $X$ is normal, an $\R$-boundary $B'$, and an $\R$-Cartier divisor $M$ on $X$ which is nef over $Z$ and such that $K_{X'}+B'+M'$ is $\R$-Cartier, where $M' : = f_*M$. We call $B'$ the boundary part and $M'$ the nef part.
\end{def*}
\begin{rem}
The notion of generalized polarized pair is a generalization of the usual notion of log pair. Indeed, we recover the latter in the case in which $X=X'$, $M=M'=0$.
\end{rem}

Now, consider a log pair $(X,B)$, and let $f:Y \rightarrow X$ be a proper birational morphism from a normal variety $Y$. Then, we can write
\begin{equation}
K_Y+f^{-1}_*B = f^*(K_X+B)+ \sum_{i=1}^n a_{E_i}(X,B) E_i,
\end{equation}
where by $f^{-1}_*B$ we denote the strict transform of $B$, the $E_i$'s are the $f$-exceptional divisors, and the coefficients $a_{E_i}(X,B)$ are the {\it discrepancies} of the $E_i$'s with respect to the pair $(X,B)$. We say that a pair is {\it log canonical}, in short {\it lc}, (respectively {\it Kawamata log terminal}, in short {\it klt}) if $a_E(X,B) \geq -1$ (respectively $a_E(X,B) > -1$) for any exceptional divisor over $X$.

In \cite{BZ}, Birkar and Zhang introduce a generalization of these measures of singularities. Consider a generalized polarized pair $(X',B'+M')$, and let $E$ be a prime divisor on any birational model of $X'$. As we are free to replace the model $X$ with a higher one, we may assume that $E$ is a prime divisor on $X$. Then, we can write
\begin{equation}
K_X + f^{-1}_*B + M = f^* (K_{X'}+B'+M') + \sum_{i=1}^n a_{E_i}(X',B'+M') E_i,
\end{equation}
where $E=E_j$ for some $1 \leq j \leq n$, and the coefficients $a_{E_i}(X',B'+M')$ are the {\it generalized discrepancies} of the $E_i$'s with respect to the generalized polarized pair $(X,B)$. We say that $(X',B'+M')$ is {\it generalized log canonical}, in short {\it glc}, (respectively {\it generalized Kawamata log terminal}, in short {\it gklt}) if $a_E(X',B'+M') \geq -1$ (respectively $a_E(X',B'+M') > -1$) for any such $E$.

In view of a possible development of a moduli theory of generalized polarized pairs, it seems reasonable to try to formulate an analog of the notion of {\it semi-log canonical} singularities (see \cite[Definition 3.13]{HK10}) in the context of generalized polarized pairs. Consider a demi-normal variety $X'$, together with a $\R$-divisors $B'$ and $M'$ such that $K_{X'} + B' + M'$ is $\R$-Cartier. Further assume that $B'$ is effective and its support is not contained in the singular locus of $X'$. Then, we say that $(X',B'+M')$ is {\it generalized semi-log canonical}, in short {\it gslc}, if $(\bar{X}',\bar{B}'+\bar{M}')$ is generalized log canonical, where $\bar{X}'$ denotes the normalization of $X'$, and $K_{\bar{X}'} + \bar{B}'+\bar{M}'$ is the pullback of $K_{X'} + B' + M'$ to $\bar{X}'$.

Since we are working over the field of complex numbers  $\C$, Hironaka's resolution of singularities holds \cite[Theorem 0.2]{KM}. In particular, given a pair $(X,B)$, we can take a {\it log resolution} $f:Y \rightarrow X$ of it. This is a proper birational morphism from a smooth variety $Y$ such that $f_*^{-1}B + \mathrm{Ex}(f)$ is a divisor with simple normal crossing support. In particular, for our purposes it will be useful to arrange families of simple normal crossing pairs. This is made precise by the following definition.

\begin{def*}
Let $X \rightarrow T$ be a morphism of varieties, and let $(X,B)$ be a pair. Then, we say that $(X,B)$ is {\it log smooth} over $T$ if $X$ is smooth, $B$ has simple normal crossing support, and every stratum of $(X,B)$ (including $X$) is smooth over $T$.
\end{def*}

\subsection{Facts about the Volume of a Generalized Polarized Pair}

Before proving our first results about the volume of generalized polarized pairs, we need one more definition.

\begin{def*}
Let $X$ be a normal variety, and consider the set of all proper birational morphisms $\pi: X_\pi \rightarrow X$, where $X_\pi$ is normal. This is a partially ordered set, where $\pi' \geq \pi$ if $\pi'$ factors through $\pi$. We define the space of {\it Weil b-divisors} as the inverse limit
\begin{equation}
\mathbf{Div}(X)\coloneqq \varprojlim_\pi \mathrm{Div}(X_\pi),
\end{equation}
where $\mathrm{Div}(X_\pi)$ denotes the space of Weil divisors on $X_\pi$. Then, we define the space of {\it $\R$-Weil b-divisors} as $\mathbf{Div}(X)_\R \coloneqq \mathbf{Div}(X)\otimes_\Z \R$.
\end{def*}

For a more detailed discussion of b-divisors in this setting, see \cite{HMX}. Here, we will just recall one construction that will appear subsequently. Consider a log pair $(X,B)$. Since a b-divisor $\mathbf{D}$ is determined by its corresponding traces $\mathbf{D}_Y$ on each birational model $Y \rightarrow X$, we can define b-divisors $\mathbf{M}_B$ and $\mathbf{L}_B$ as follows. For every birational morphism $\pi : Y \rightarrow X$, define $B_Y$ by $K_Y+B_Y : =\pi^*(K_X+B)$. Then, we declare
\begin{equation}
\mathbf{M}_{B,Y} \coloneqq \pi_*^{-1}B+\mathrm{Ex}(\pi), \quad \mathbf{L}_{B,Y} \coloneqq B_Y \vee 0.
\end{equation}
Here $\vee$ denotes the maximum between two divisors; similarly, we will use $\wedge$ to denote the minimum between two divisors. The b-divisor $\mathbf{M}_B$ encodes the information of the strict transforms and the exceptional divisors on all higher models, while $\mathbf{L}_B$ records the effective part of the log pullback of $B$ on any higher model.

Now, we are ready to state the following lemma, which is a generalization of \cite[Lemma 2.2.1]{HMX}.

\begin{lemma} \label{volumes lemma}
Let $f:X \rightarrow W$ and $g: Y\rightarrow X$ be birational morphisms of normal projective varieties, $D$ an $\mathbb{R}$-divisor on $X$, and $M$ an $\mathbb{R}$-Cartier $\mathbb{R}$-divisor on $W$. Then the following statements hold.
\begin{itemize}
\item[(1)] $\vol ( W,f_* D) \geq \vol (X,D)$.
\item[(2)] If $D$ is $\mathbb{R}$-Cartier and $G$ is an $\mathbb{R}$-divisor on $Y$ such that $G-g^*D \geq 0$ is effective and $g$-exceptional, then $\vol(Y,G)=\vol(X,D)$. In particular, if $(X,B)$ is a projective log canonical pair and $f:Y \rightarrow X$ a birational morphism, then
\begin{equation}
\begin{split}
\vol (X,K_X+B+N)&=\vol(Y,K_Y+\mathbf{L}_{B,Y}+g^*N)\\&=\vol(Y,K_Y+\mathbf{M}_{B,Y}+g^*N),
\end{split}
\end{equation}
where $N$ is any $\mathbb{R}$-Cartier $\mathbb{R}$-divisor on $X$.
\item[(3)] If $D \geq 0$, $(W,f_* D)$ has simple normal crossings, and $L= \mathbf{L}_{f_*D,X}$, then
\begin{equation}
\vol (X,K_X+D+f^*M)=\vol(X,K_X+D\wedge L +f^*M).
\end{equation}
\item[(4)] If $(X,B)$ is a log canonical pair and $X \dashrightarrow X'$ is a birational contraction of normal projective varieties, then
\begin{equation}
\vol (X',K_{X'}+\mathbf{M}_{B,X'})\geq \vol (X,K_X+B).
\end{equation}
If, moreover, $f : X \rightarrow W$ and $h : X' \rightarrow W$ are morphisms, then we have
\begin{equation}
\vol (X',K_{X'}+\mathbf{M}_{B,X'}+h^*M)\geq \vol (X,K_{X}+B+f^*M).
\end{equation}
In addition, if the centre of every divisor in the support of $B\wedge \mathbf{L}_{f_*B,X}$ is a divisor on $X'$, and $(W,f_*B)$ has simple normal crossings, then we have equality
\begin{equation}
\vol (X',K_{X'}+\mathbf{M}_{B,X'}+h^*M)=\vol (X,K_{X}+B+f^*M).
\end{equation}
\end{itemize}
\end{lemma}
\begin{rem}
In part (4), $\mathbf{M}_{B,X'}$ is constructed as follows: We take a common resolution $X''$, dominating both $X$ and $X'$. Then, $\mathbf{M}_{B,X'}$ is the pushforward of $\mathbf{M}_{B,X''}$ to $X'$.
\end{rem}
\begin{dimo1*}
As for (1) and (2), see \cite[Lemma 2.2.1]{HMX}. As for (3) and (4), the proof is a slight generalization of what is in \cite[Lemma 2.2.1]{HMX}.

Now we will prove (3). One has the obvious inclusion
\begin{equation}
H^0(X,\Ox_X(m(K_X+D+f^*M)))\supset H^0(X,\Ox_X(m(K_X+D\wedge L+f^*M)))
\end{equation}
for all $m \geq 0$, which leads to the inequality
\begin{equation}
\vol (X,K_X+D+f^*M)\geq \vol(X,K_X+D\wedge L +f^*M). 
\end{equation}
Now, we have
\begin{equation}
K_X+L+f^*M=f^*(K_W+f_*D+M)+E,
\end{equation}
where $L \wedge E =0$, $L \geq 0$ and $E \geq 0$. Furthermore, $E$ is $f$-exceptional. This guarantees
\begin{equation}
\begin{split}
H^0(X,\Ox_X(m(K_X+D+f^*M))) &\subset f^*H^0(W,\Ox_W(m(K_W+f_*D+M))) \\
&=H^0(X,\Ox_X(m(K_X+L+f^*M))),
\end{split}
\end{equation}
where the first inclusion follows from (1), and the second equality from (2). Thus, sections of $H^0(X,\Ox_X(m(K_X+D+f^*M)))$ vanish along $D-D\wedge L$, so the claim follows.

Now, we are left with proving (4). The first part of the statement is proved in \cite[Lemma 2.2.1]{HMX}. The second part follows from the first one, and part (2). As for the third part of the statement, set $B' \coloneqq \mathbf{M}_{B,X'}$, and let $X''$ be a resolution of indeterminacies of $X \dashrightarrow X'$. Then, the proof of \cite[Lemma 2.2.1]{HMX} provides us with
\begin{equation} \label{eq to refer to}
\mathbf{M}_{B',X''}\wedge\mathbf{L}_{h_*B',X''}=\mathbf{M}_{B,X''}\wedge\mathbf{L}_{f_*B,X''}.
\end{equation}
Therefore, we have
\begin{equation}
\begin{split}
\vol (X,K_X+B+f^*M)&= \vol (X'',K_{X''}+\mathbf{M}_{B,X''}\wedge \mathbf{L}_{f_*B,X''}+f^*M)\\
&= \vol (X'',K_{X''}+\mathbf{M}_{B',X''}\wedge \mathbf{L}_{h_*B',X''}+f^*M)\\
&=\vol (X',K_{X'}+B'+h^*M),
\end{split}
\end{equation}
where the first equality follows from (2) and (3) together, the second one from identity (\ref{eq to refer to}), and the third one from (2) and (3) again. \hfill $\square$
\end{dimo1*}

Now, we can deal with the volume of birational models living over a fixed simple normal crossing pair. In particular, we are ready to provide the proof of Theorem \ref{DCCcase1}, which is a generalization of \cite[Theorem 3.1.2]{HMX}.

\begin{dimo1*}[{Theorem \ref{DCCcase1}}]
Suppose, by contradiction, that there is an infinite sequence $\lbrace (X_i,B_i) \rbrace_{i \geq 1} \subset \mathfrak{D}$ such that $\vol(X_i,K_{X_i}+B_i+f_i^*M)$ is strictly decreasing. By Lemma \ref{volumes lemma}, we know that
\begin{equation}
\vol(X_i,K_{X_i}+B_i+f_i^*M) \leq \vol(Z,K_{Z}+f_{i,*}B_i+M),
\end{equation}
and in case of strict inequality, there is some $E \subset \mathrm{Ex}(f_i)$ such that
\begin{equation}
\mult_E (K_{X_i}+B_i+f_i^*M) < \mult_E (f_i^*(K_{Z}+f_{i,*}B_i+M)),
\end{equation}
which is equivalent to
\begin{equation}
\mult_E (K_{X_i}+B_i) < \mult_E (f_i^*(K_{Z}+f_{i,*}B_i)).
\end{equation}
By Lemma \ref{volumes lemma}, which now takes care of generalized polarized pairs involving $M$, following the lines of \cite[Theorem 3.1.2]{HMX}, we may assume that all the $f_i$'s are toroidal. Furthermore, after passing to a subsequence, we have that the $\mathbf{M}_{B_i}$'s form a non-decreasing sequence of b-divisors converging to $\mathbf{B} \coloneqq \lim \mathbf{M}_{B_i}$.

By \cite[Lemma 3.1.1]{HMX}, there is a reduction $(\mathbf{B'},Z')$ of $(\mathbf{B},Z)$ such that $\mathbf{B'}\geq \mathbf{L}_{\mathbf{B'}_{Z'}}$. 
Furthermore, by construction of a reduction \cite[pp. 22-23]{HMX} and part (3) of Lemma \ref{volumes lemma}, we have the equality
\begin{equation}
\vol(Z',K_{Z'}+\mathbf{B'}_{Z'}+f^*M)=\vol(Z',K_{Z'}+\mathbf{B}_{Z'}+f^*M).
\end{equation}

Since the boundaries $B_i$ in the varieties $X_i$ are log canonical, up to replacing $(X_i,B_i)$ with a higher model, we may assume the $f_i$'s factor through $Z'$.

For $i$ sufficiently large, we may assume $\mathbf{M}_{B_i,Z'}\leq \mathbf{B}_{Z'}$, since the coefficients are in a DCC set and the support is fixed on the model $Z'$. Thus, we have
\begin{equation}
\begin{split}
\vol(X_i,K_{X_i}+B_i+f_i^*M) &\leq \vol(Z',K_{Z'}+\mathbf{M}_{B_i,Z'}+f^*M)\\
&\leq \vol(Z',K_{Z'}+\mathbf{B}_{Z'}+f^*M)\\
&= \vol(Z',K_{Z'}+\mathbf{B'}_{Z'}+f^*M),
\end{split}
\end{equation}
where the first inequality follows from part (1) of Lemma \ref{volumes lemma}.

Since the volumes $\vol(X_i,K_{X_i}+B_i+f_i^*M)$ form a strictly decreasing sequence, for any $i >>1$ there is $\epsilon > 0$ such that
\begin{equation}
\begin{split}
\vol(X_j,K_{X_j}+B_i+f_j^*M)&< \vol(X_i,K_{X_i}+B_i+f_i^*M) \\ &\leq \vol(Z',K_{Z'}+(1-\epsilon)\mathbf{B'}_{Z'}+f^*M)
\end{split}
\end{equation}
for any $j>i$. Now, let $Z'' \rightarrow Z$ be a toroidal morphism that extracts all divisors $E$ with $a_E(Z',(1-\epsilon )\mathbf{B'}_{Z'})<0$. Then we have
\begin{equation}
\mathbf{M}_{B_j,Z''}\geq \mathbf{L}_{(1-\epsilon)\mathbf{B'}_{Z'},Z''}
\end{equation}
for $j>>0$, which implies
\begin{equation}
\mathbf{M}_{B_j,X_j}\geq \mathbf{L}_{(1-\epsilon)\mathbf{B'}_{Z'},X_j}.
\end{equation}
Therefore, we get
\begin{equation}
\begin{split}
\vol(X_j,K_{X_j}+B_j+f_j^*M) &\geq \vol(X_j,K_{X_j}+\mathbf{L}_{(1-\epsilon)\mathbf{B'}_{Z'},X_j}+f_j^*M)\\
&= \vol(Z'',K_{Z''}+\mathbf{L}_{(1-\epsilon)\mathbf{B'}_{Z'},Z''}+g^*M) \\
&= \vol(Z',K_{Z'}+(1-\epsilon)\mathbf{B'}_{Z'}+f^*M) ,
\end{split}
\end{equation}
where $g:Z'' \rightarrow Z$. This provides the required contradiction. \hfill $\square$
\end{dimo1*}

Now, we can proceed with the proof of the second major general statement, namely the deformation invariance of volumes in the case of generalized polarized pairs.

\begin{dimo1*}[{Theorem \ref{DefInvariance}}]
We need to compare the volumes pairwise. Therefore, after base change, we may assume that $T$ is a smooth curve containing $x_i$ and $x_j$, for some fixed $i,j \in \N$. Now, let $\mathcal{H}$ be an ample divisor on $\mathcal{X}$. Then, $\mathcal{A}_n \coloneqq \mathcal{M}+\frac{1}{n}\mathcal{H}$ is ample over $x_i$ and $x_j$. Thus, by openness of such a property in families, we can find $x_i,x_j \in T_n \subset T$ open such that $\mathcal{A}_n$ is relatively ample over $T_n$.

Hence, up to the pullback of a sufficiently ample class on $T_n$, the divisor $\mathcal{A}_n$ can be considered ample. After shrinking $T_n$, we may assume $\mathcal{A}_n$ is ample itself.

If we choose $\mathcal{D}_n \equiv \mathcal{A}_n$ generically enough, we may assume that $\mathcal{B}+\mathcal{D}_n$ has simple normal crossings, and, over an open subset $U_n \subset T_n$ containing $x_i$ and $x_j$, $(\mathcal{X},\mathcal{B}+\mathcal{D}_n)$ is log smooth. Therefore, by \cite[Theorem 2.6.2]{HMX}, we have invariance of plurigenera over $U_n$. In particular, as $n$ varies, we have
\begin{equation}
\vol(X_i,K_{X_i}+B_i+\mathcal{A}_n|_{X_i})=\vol(X_j,K_{X_j}+B_j+\mathcal{A}_n|_{X_j}).
\end{equation}
As $n \to +\infty$, we have that $\mathcal{A}_n \to \mathcal{M}$ in the Euclidean topology of $N^1(\mathcal{X})_\R$. This implies $\mathcal{A}_n|_{{X}_i} \to M_i$ and $\mathcal{A}_n|_{{X}_j} \to M_j$ in the Euclidean topologies of $N^1({X}_i)_\R$ and $N^1({X}_j)_\R$, respectively. Therefore, by continuity of the volume, we have
\begin{equation}
\vol(X_i,K_{X_i}+B_i+M_i)=\vol(X_j,K_{X_j}+B_j+M_j).
\end{equation}
Repeating this argument for all pairs $i,j \in \N$ leads to the conclusion. \hfill $\square$
\end{dimo1*}

\subsection{Some facts towards the DCC property}

In this paragraph, we collect some technical statements that are used to prove Theorem \ref{MainResult}, yet are already phrased in the more general setting of arbitrary dimension.

Now, we are interested in showing that a particular set of volumes satisfies the DCC property, and we have a tool that allows us to compare volumes in families, namely Theorem \ref{DefInvariance}. Therefore, the natural strategy is to put the log pairs of interest in a family, which will be made more precise with the following definition.

\begin{def*}
Let $\mathfrak{D}$ be a set of log pairs. We say that $\mathfrak{D}$ is {\it bounded} (respectively {\it log birationally bounded}) if there is a log pair $(\mathcal{X},\mathcal{B})$ with $\mathcal{B}$ reduced, and there is a projective morphism $\mathcal{X} \rightarrow T$, where $T$ is of finite type, such that for every log pair $(X,B) \in \mathfrak{D}$ there are a closed point $t \in T$  and a morphism $f:\mathcal{X}_t \rightarrow X$ inducing an isomorphism $(X,B_{\mathrm{red}}) \cong (\mathcal{X}_t,\mathcal{B}_t)$ (respectively, a birational map such that the support of $\mathcal{B}_t$ contains the support of the strict transform of $B$ and of any $f$-exceptional divisor).
\end{def*}

\begin{rem} \label{BoundLogBir}
By a standard Hilbert scheme argument, a set of log pairs $\mathfrak{D}$ is bounded if for any $(X,B) \in \mathfrak{D}$ there is a very ample divisor $H$ such that $H^{\dim X} \leq C$ and $B_\red \cdot H^{\dim X - 1}\leq C$ for some constant $C$ \cite[Remark 2.7.3]{HMX}. Furthermore, this provides us with a strategy to show that a set of log pairs is log birationally bounded. Indeed, assume that for any log pair $(X,B) \in \mathfrak{D}$ there is a birational morphism $f:X \rightarrow Y$ induced by a free linear series $|H|$. Then, bounding the intersection products $H^{\dim X}$ and $B_\red \cdot H^{\dim X -1}$ is equivalent to bounding $A^{\dim Y}$ and $f_*B_\red \cdot A^{\dim Y -1}$, where $H=f^*A$ and $A$ is very ample on $Y$. Bounding the latter products is equivalent to saying that the pairs $(Y,f_*B)$ are bounded; since there are no exceptional divisors for $Y \dashrightarrow X$, the last condition implies that $\mathfrak{D}$ is log birationally bounded.
\end{rem}

On the other hand, we will see that, in the setting of generalized polarized pairs, the concept of boundedness is more complicated, since the best we can hope for is to have a divisor $\mathcal{M}$ on the bounding family $\mathcal{X}$ such that $\mathcal{M}_t \sim_\Q M_t$ (or even $\mathcal{M}_t \equiv M_t$), where $M_t$ is the nef class living on $\mathcal{X}_t$.

In view of Theorem \ref{BirkarZhangBirational}, we have natural candidates for the maps $f$ that will lead to the construction of a log birationally bounded family. Since we are interested in volumes of log canonical pairs, we can always resolve the indeterminacies of the map going to a higher model, and create a new pair having the same volume. We make this precise in the next statement.

\begin{lemma} \label{smoothModel}
Assume that $(X,B+M)$ is a generalized polarized pair such that $M$ is $\R$-Cartier, and $X=X'$, $M=M'$, and $(X,B)$ is log canonical. Furthermore, assume that $|\lfloor m_0(K_X+B+M) \rfloor|$ defines a birational map. Then, there exist a smooth birational model $f:Y \rightarrow X$ and a boundary $D$ such that $(Y,D)$ is log canonical, $|\lfloor m_0(K_Y+D+f^*M) \rfloor|$ defines a birational morphism, and
\begin{equation}
\vol (X,K_X+B+M) = \vol (Y,K_Y+D+f^*M).
\end{equation}
Furthermore, we have that $\coeff (D) \subset \coeff(B) \cup \lbrace 1 \rbrace$.
\end{lemma}
\begin{rem}
Since the nef part in the higher model is just the pullback $f^*M$ of $M$, in case $M$ is $\Q$-Cartier, the Cartier index of the nef part remains unchanged.
\end{rem}
\begin{dimo1*}
Note that, since $K_X+B+M$ is $\R$-Cartier, and $M$ is $\R$-Cartier, then $K_X+B$ is $\R$-Cartier. Thus, it makes sense for $(X,B)$ to be log canonical. Now, let $Y$ be a log resolution of both $(X,B)$ and $|\lfloor m_0(K_X+B+M) \rfloor|$. Let $D \coloneqq \mathbf{M}_{B,Y}$. Then, by Lemma \ref{volumes lemma}, we have
\begin{equation}
\vol (X,K_X+B+M) = \vol (Y,K_Y+D+f^*M).
\end{equation}
Now, since $f^*|\lfloor m_0(K_X+B+M) \rfloor|$ defines a birational morphism, and $f^*|\lfloor m_0(K_X+B+M) \rfloor| = |\lfloor m_0(K_Y+D+f^*M) \rfloor|$, we are done. \hfill $\square$
\end{dimo1*}

Now, we can make a first step towards the log birational boundedness of the pairs involved in the statement of Theorem \ref{MainResult}. In particular, the following shows how to bound the boundary part $B$.

\begin{prop} \label{boundHB}
Let $(X,B)$ be a log canonical pair with $X$ smooth, and $M$ a nef $\R$-Cartier $\R$-divisor on $X$. Assume that $|\lfloor m_0( K_X +B+ M) \rfloor|=|H|+F$ defines a birational morphism. Furthermore, assume that $B \in \Lambda$, where $\Lambda \subset [0,1]$ and $\Lambda \cap (0,1]$ has a positive minimum $\delta$. Then, we have
\begin{equation}
B_\red \cdot H^{n-1} \leq \left(2 \left( \frac{m_0(n+1)+2}{\delta}+2(2n+1)m_0+2 \right) \right)^n \vol (X,K_X+B+M),
\end{equation}
where $n=\dim X$.
\end{prop}

\begin{dimo1*}
By \cite[Lemma 2.7.6]{HMX}, we have
\begin{equation}
|K_X+(n+1)\lfloor m_0(K_X+M+B) \rfloor|\neq \emptyset.
\end{equation}
Also, we have that $\delta B_\red \leq B$. Furthermore, since $M$ is nef and $K_X+M+B$ is big, the divisor $K_X+2M+B$ is big, and then we have
\begin{equation} \label{Kplus2MplusBeffective}
K_X+2M+B \sim_\R E \geq 0
\end{equation}
for some effective $\R$-divisor $E$. Therefore, there exists an effective $\R$-divisor $C \geq 0$ such that
\begin{equation} \label{BplusEffective}
B_\red + C \sim_\R \frac{m_0(n+1)+2}{\delta}(K_X+M+B).
\end{equation}
Therefore, we have the following chain of inequalities
\begin{equation}
\begin{split}
B_\red \cdot (2(2n+1)H)^{n-1}  \leq 2^n \vol (X,K_X+B_\red+2(2n+1)H)\\
 \leq 2^n \vol \left(X,K_X+\left(\frac{m_0(n+1)+2}{\delta}+  2(2n+1)m_0\right)(K_X+M+B)\right)\\
 \leq 2^n \vol \left(X,K_X+B+ \left(\frac{m_0(n+1)+2}{\delta}  +2(2n+1)m_0\right)(K_X+M+B)\right)\\
 \leq 2^n \vol \left(X,2K_X+2M+2B+\left(\frac{m_0(n+1)+2}{\delta}+2(2n+1)m_0\right)(K_X+M+B)\right)\\
 = 2^n \vol \left(X,\left(\frac{m_0(n+1)+2}{\delta}+2(2n+1)m_0+2\right)(K_X+M+B)\right)\\
 = \left(2\left(\frac{m_0(n+1)+2}{\delta}+2(2n+1)m_0+2\right)\right)^n \vol (X,(K_X+M+B)),
\end{split}
\end{equation}
where the first inequality follows from \cite[Lemma 2.7.5]{HMX}, the second one from equation \ref{BplusEffective}, and the third one from equation \ref{Kplus2MplusBeffective}. \hfill $\square$
\end{dimo1*}

\section{Some Facts about Surfaces}

As previously mentioned, the morphisms coming from Theorem \ref{BirkarZhangBirational} will be the ones to provide the log birational boundedness of the log pairs of interest. On the other hand, there are several subtleties to be aware of in dealing with the nef part $M$.

Consider a birational morphism $f:X \rightarrow Y$ of normal projective varieties, and let $M$ be a nef $\R$-Cartier divisor on $X$. As a first problem, there is no guarantee that $f_*M$ is $\R$-Cartier. More importantly, even assuming that $f_*M$ is $\R$-Cartier, we do not know whether $f_*M$ is nef as well. Assuming the latter condition, the best we can say is that, by the Negativity Lemma \cite[Lemma 3.39]{KM}, $f^*f_*M-M$ is effective and exceptional.

From this perspective, there are a few reasons why there is a more explicit way to tackle our problem in the case of surfaces. First of all, most of the time, we can work in the smooth category, where birational morphisms have a very explicit characterization \cite[Theorem II.11]{BEAU}. Furthermore, even in the case that the target is not $\Q$-factorial, we can still make sense of the pullback of Weil divisors, since the intersection matrix of the exceptional curves is negative definite \cite[Lemma 3.40]{KM}. Also, as we will recall in Lemma \ref{push nef}, if $X$ and $Y$ are smooth surfaces and $M$ is nef, then $f_* M$ is nef as well. Last but not least, with Theorem \ref{structureSurfaceNef} we have an explicit way to determine how high over $Y$ the difference between $M$ and $f_*M$ lives.

\subsection{Towards the Proof of Theorem \ref{structureSurfaceNef}}

Now, we will recall a few facts that will be important in the proof of Theorem \ref{structureSurfaceNef}.

\begin{lemma} \label{blow up point}
Let $X$ be a smooth surface, and $M$ a nef Cartier divisor on $X$. Let $X' \coloneqq Bl_p X$, for some closed point $p \in X$. Denote by $\pi : X' \rightarrow X$ the induced morphism, and by $E$ the $\pi$-exceptional divisor. Then, $\pi^*M-E$ is nef if and only if $M \cdot C \geq \mult_p C$ for any irreducible curve $C \subset X$ through $p$.
\end{lemma}
\begin{dimo1*}
Let $D$ be an irreducible curve on $X'$. If $D \cap E = \emptyset$, then $D=\pi^* \pi_* D$. Therefore, $(\pi^*M-E)\cdot D= M \cdot \pi_* D \geq 0$. On the other hand, we have $(\pi^*M-E)\cdot E = -E^2=1$.

Now, consider an irreducible curve $C' \subset X'$ such that $C' \cap E$ is a finite set. Also, write $C \coloneqq \pi_* C'$. Then, we have $\pi^*C=C'+(\mult_p C) E$. Thus, we have
\begin{equation}
(\pi^*M-E)\cdot C'= (\pi^*M-E)\cdot (\pi^*C-(\mult_pC)E)= C \cdot M - \mult_p C.
\end{equation}
The statement follows. \hfill $\square$
\end{dimo1*}

As an immediate consequence of Lemma \ref{blow up point}, we get the following statement.

\begin{cor} \label{blow up case}
Let $X$ be a smooth surface, and $M$ be a nef Cartier divisor. Assume there is an irreducible curve $C$ such that $M \cdot C=0$. Then, for any $p \in C$, $\pi_p^*M-E_p$ is not nef, where $\pi_p : X_p \rightarrow X$ is the blow-up at $p$ of $X$ and $E_p$ is the $\pi_p$-exceptional divisor.
\end{cor}

Now, we will start to compare the volume of nef divisors living on different models.

\begin{lemma} \label{drop volume}
Let $\pi : X_1 \rightarrow X$ be a birational morphism of smooth surfaces. Assume there are nef Cartier divisors $M$ on $X$ and $M_1$ on $X_1$ such that $\pi_*M_1=M$. Then, we have $\pi^*M=M_1+E$, where $E$ is effective and $\pi$-exceptional. Furthermore, $M^2=M_1^2$ if and only if $E=0$.
\end{lemma}

\begin{dimo1*} Since $\pi_*M_1=M$, we have that $\pi^*M=M_1+E$, where $E$ is a $\pi$-exceptional divisor. Since $-E\sim_{\pi} M_1$ is $\pi$-nef, by the Negativity Lemma \cite[Lemma 3.39]{KM} $E$ is effective. Now, we know that the intersection matrix of the $\pi$-exceptional curves is negative definite \cite[Lemma 3.40]{KM}. Therefore, $E^2 \leq 0$, and $E^2=0$ if and only if $E=0$. Furthermore, $M_1^2=(\pi^*M-E)^2=M^2+E^2$. This, together with the analysis of $E^2$, concludes the proof. \hfill $\square$
\end{dimo1*}

Next, we recall the well known fact that the pushforward of a nef divisor remains nef in the case of smooth surfaces.

\begin{lemma} \label{push nef}

Let $f:Y \rightarrow X$ be a birational morphism of smooth surfaces, and let $M$ be a nef divisor on $Y$. Then, $f_*M$ is nef on $X$.

\end{lemma}

\begin{dimo1*}
Since $f:Y \rightarrow X$ factors as finitely many blow-downs of $(-1)$-curves \cite[Theorem II.11]{BEAU}, we may assume $Y=Bl_pX$ for some $p \in X$. Now, let $E$ be the $f$-exceptional divisor. Then, $E^2=-1$. Furthermore, since $M$ is nef, by the Negativity Lemma \cite[Lemma 3.39]{KM}, $f^*f_*M=M+aE$ for some $a \geq 0$.

Now, let $C \subset X$ be an irreducible curve. Then, $f^*C=C'+(\mult_p C)E$, where $C'$ is the strict transform of $C$. Also, we have that $C'\cdot E=\mult_p C$. Therefore, we have
\begin{equation}
\begin{split}
C \cdot f_*M&=f^*C \cdot f^* f_* M\\&=( C'+(\mult_p C)E)\cdot(M+aE)\\&=C'\cdot M+(\mult_pC)E\cdot M\\& \geq 0.
\end{split}
\end{equation}
This concludes the proof. \hfill $\square$
\end{dimo1*}

\begin{oss} \label{intersection goes up}
The proof also shows the following facts. First, if $C$ is disjoint from the centers of the blow-ups, then $C \cdot f_*M=C'\cdot M$. Furthermore, consider the setting of one simple blow-up $\pi_p : Bl_pX \rightarrow X$, and assume that $M \cdot E >0$. Then, for $C \subset X$ through $p$, we have that $C\cdot \pi_{p,*}M > C'\cdot M$.
\end{oss}

\begin{rem} \label{push nef singular}
The statement of Lemma \ref{push nef} goes through in the case the varieties $X$ and $Y$ are just normal, assuming that $f_*M$ is $\Q$-Cartier. By continuity, we may assume that $M$ is ample. Then, a multiple of $f_*M$ is semiample off a finite set of points. Therefore, by a Theorem of Zariski \cite[p. 132]{LAZ1}, such a multiple is honestly semiample, and hence nef.
\end{rem}

Now, we can combine Corollary \ref{blow up case} with Lemma \ref{drop volume} to get the key tool for the proof of Theorem \ref{structureSurfaceNef}.
\begin{prop} \label{trivial blow up}
Assume we have a birational morphism between smooth surfaces $f:X' \rightarrow X$, and let $M'$ and $M$ be nef divisors on $X'$ and $X$ respectively such that $f_*M'=M$. Also, assume there is an irreducible curve $C \subset X$ such that $M \cdot C=0$. Then $M' \cdot E=0$ for any $f$-exceptional curve such that $f(E) \in C$. 
\end{prop}
\begin{dimo1*}
By the structure theorem for birational morphisms between smooth surfaces \cite[Theorem II.11]{BEAU}, we have that $X'$ is obtained from $X$ by a finite sequence of blow-ups. Also, blow-ups of distinct points (i.e. not infinitely close) on $X$ do not interfere with each other: We can blow them up independently.  Therefore, by induction on the number of centers of the blow-up, we may assume that $X'$ is obtained blowing up one closed point $p \in X$ and points infinitely close to it.

Suppose that we blow up in order the points $p_1 \coloneqq p, p_2,  \ldots, p_n$. We inductively define surfaces $X_i$ and morphisms $\pi_i$ as follows: Let $X_0 \coloneqq X$, and denote by $\pi_1 : X_1 \rightarrow X_0$ the blow-up of $X_0$ at $p_1$. Then, $X_{i+1}$ and $\pi_{i+1}$ are defined by the blow-up $\pi_{i+1}:X_{i+1}\rightarrow X_i$ of $X_i$ at $p_i \in X_i$. In particular, we have $X_n=X'$. Also, denote by $f_i: X' \rightarrow X_i$ and $g_i : X_i \rightarrow X$ the morphisms naturally induced by this construction. Finally, define $M_i \coloneqq f_{i,*}M'$. In particular, we have $M_0=M$. The picture is summarized as follows.

\begin{center}
\begin{tikzpicture}
  \matrix (m) [matrix of math nodes,row sep=3em,column sep=4em,minimum width=2em]
  {
     (X'=X_n,M'=M_n) & (X_i,M_i) \\
      & (X_{i-1},M_{i-1}) \\
      & (X=X_0,M=M_0) \\};
  \path[-stealth]
    (m-1-1) edge node [above] {$f_i$} (m-1-2)
            edge node [left] {$f=f_0$} (m-3-2)
    (m-1-2) edge node [right] {$\pi_i$} (m-2-2)
    (m-2-2) edge node [right] {$g_{i-1}$} (m-3-2);
\end{tikzpicture}
\end{center}

Now, assume $p \in C$, where $C$ is an irreducible curve, and $M \cdot C = 0$. We prove the statement by induction on the number of blow-ups $n$. More precisely, after the above reductions, we will prove that $M'=f^*M$. If $n=1$, this follows by the proof of Lemma \ref{push nef}. Notice that in this case $M_1=g_1^*M$ holds. Then, assume this holds for $m<n$, and we will show it holds for $n$ as well. In particular, the inductive hypothesis applies to $(X_{n-1},M_{n-1}) \rightarrow (X,M)$. Therefore, we have $M_{n-1} = g_{n-1}^*M$, which implies $M_{n-1} \cdot E =0$ for any $g_{n-1}$-exceptional curve. Since $p_n \in E$ for some $g_{n-1}$-exceptional curve $E$, we are in the setting of Corollary \ref{blow up case}. Therefore, $M'=f_{n-1}^* M_{n-1}=f^* M$, and this concludes the proof.   \hfill $\square$
\end{dimo1*}

Now we are ready to prove Theorem \ref{structureSurfaceNef}

\begin{dimo1*}[{Theorem \ref{structureSurfaceNef}}]
We can proceed by induction on $M^2$, the base of the induction $M^2=0$ being provided by Lemma \ref{drop volume}. Since $f$ is a birational morphism, we know that $X'$ is obtained by blowing up finitely many points on $X$, possibly infinitely close \cite[Theorem II.11]{BEAU}. Therefore, there is at least one $f$-exceptional $(-1)$-curve $E$. If $M_Y \cdot E =0$, we can blow $E$ down via $\pi: Y \rightarrow Y'$. By the Negativity Lemma \cite[Lemma 3.39]{KM} and Lemma \ref{blow up point}, we know that $M_Y=\pi^* \pi_* M_Y$. We can repeat this process just finitely many times; call $X'$ the resulting model. Thus, either $X'=X$, in which case we are done, or all the $f'$-exceptional $(-1)$-curves $E$ are such that $M_{X'} \cdot E >0$, where $M_{X'}$ denotes the pushforward of $M_Y$ to $X'$, and $f'$ is the induced morphism $f':X'\rightarrow X$.

\begin{claim}
The $X'$ above constructed satisfies the claims of the thesis. 
\end{claim}

The claim is obvious when $X'=X$. Now, in the case that $X'$ is not isomorphic to $X$, there exists a point $p \in X$ such that $f$ factors through $\pi_p: Bl_p X \rightarrow X$. Call $h$ the map $h: X' \rightarrow Bl_pX$. Also, write $M_p \coloneqq h_* M_{X'}$. By the Negativity Lemma \cite[Lemma 3.39]{KM}, $M_p = \pi_p^* M - a E_p$, where $a \geq 0$ and $E_p$ is the $\pi_p$-exceptional divisor. If $a > 0$, then $M_p^2 \leq M^2 - 1$. In such a case, the inductive hypothesis applies to $h:(X',M') \rightarrow (Bl_pX,M_p)$.

Now, we are going to show that $a=0$ is not admissible. Call $C_p \coloneqq h_*^{-1}E_p$, and notice that, in this setting, $M_p=\pi_p^* M$. Also, we have that $M'=h^*M_p-F$, where $F \geq 0$ is $h$-exceptional\footnote{Notice that the assumption $a=0$ together with $M' \cdot E >0$ for all $f'$-exceptional $(-1)$-curves implies that $h$ is not an isomorphism.}. In particular, the supports of $C_p$ and of $F$ have no common components.

First, assume $C_p$ is a $(-1)$-curve. Then, $h^{-1}$ is an isomorphism in a neighborhood of $E_p$. Therefore, the supports of $C_p$ and $F$ are disjoint. Thus, we have
\begin{equation}
M' \cdot C_p= (h^*M_p-F)\cdot C_p = h^*M_p \cdot C_p = f^*M \cdot C_p=0.
\end{equation}
On the other hand, by construction, $M'$ is positive against all the $f'$-exceptional $(-1)$-curves. In particular, we must have $M' \cdot C_p >0$. Therefore, this setting is not possible.

The last case we have to rule out is the following: There exists an $h$-exceptional $(-1)$-curve $L$ such that $h(L) \in E_p$. Since we are assuming $a=0$, by Proposition \ref{trivial blow up}, we have that $M'$ is trivial against any curve mapped to $E_p$. This leads to a contradiction. Therefore, the case $a=0$ is not admissible, and this concludes the proof. \hfill $\square$
\end{dimo1*}

\subsection{Some Bounds in the Surface Case}

Here we include a couple of results that go in the direction of bounding suitable representatives of the pairs in the statement of Theorem \ref{MainResult}. First, we consider the intersection product between the nef part $M$ and a general element of the free part of the linear series defining a morphism as in Theorem \ref{BirkarZhangBirational}.

\begin{prop} \label{SurfaceBoundHM}
Let $X$ be a smooth surface, $B \geq 0$ a boundary divisor, and $M$ a nef Cartier divisor. Assume there is $m_0 \in \N$ such that $|\lfloor m_0(K_X+B+M) \rfloor|=|H|+F$ defines a birational morphism. 
Then, we have
\begin{equation}
H\cdot (M+3H) \leq \frac{(6m_0+1)^2}{2} \vol(X, K_X+B+M).
\end{equation}
\end{prop}
\begin{dimo1*}
Since $|H|$ is basepoint-free, its general element is smooth; also, we may assume that such an element is not contained in any given finite set of curves. This being said, assume that $C \in |H|$ is such a general element. Also, since $M$ is nef and $H$ is nef and big, we have that $M+kH$ is nef and big for any $k>0$.

Now, let $k \in \N$ and $N$ be any nef Cartier divisor. We consider the following short exact sequence
\begin{equation}
0 \rightarrow \Ox_X(K_X+N+kH) \xrightarrow{\cdot C} \Ox_X(K_X+N+kH +H) \rightarrow \Ox_C(K_C+N_C+kH_C) \rightarrow 0,
\end{equation}
where $N_C$ and $H_C$ denote the restrictions of $N$ and $H$ to $C$. Since $H^2\geq 1$, and $N \cdot H \geq 0$, then
\begin{equation}
\deg (K_C+N_C+kH_C)\geq 2g(C)-2+k.
\end{equation}
Therefore, by Riemann-Roch, $\Ox_C(K_C+N_C+kH_C)$ is free for $k\geq 2$, and it is very ample for $k \geq 3$. Furthermore, by Kawamata-Viehweg vanishing, we have $H^1(X,\Ox_X(K_X+N+kH))=0$ for $k\geq 1$. Therefore, $H^0(X,\Ox_X(K_X+N+kH))\neq 0$ for $k \geq 2$.  In particular, letting $N=0$, we have $|K_X+3H|\neq \emptyset$.

Now, pick an element $D\in |(m-1)(K_X+3H)|$; we may assume that the supports of $D$ and $C$ share no components. Then, multiplication by $D$ leads to the following commutative diagram
\begin{center}
\begin{tikzpicture}
  \matrix (m) [matrix of math nodes,row sep=3em,column sep=4em,minimum width=2em]
  {
     \Ox_X(K_X+m(M+3H)+H) & \Ox_C(K_C+m(M_C+3H_C)) \\
     \Ox_X(m(K_X+M+6H)-2H) & \Ox_C(m(K_C+M_C+5H_C)-2H_C) \\};
  \path[-stealth]
    (m-1-1) edge node [left] {$\cdot D$} (m-2-1)
            edge node [below] {} (m-1-2)
    (m-2-1) edge node {} (m-2-2)
    (m-1-2) edge node [right] {$\cdot D|_C$} (m-2-2);
\end{tikzpicture}
\end{center}
where $m \in \N$. By construction, the induced maps 
\begin{equation}
H^0(C,\Ox_C(K_C+m(M_C+3H_C))) \rightarrow H^0(C,\Ox_C(m(K_C+M_C+5H_C)-2H_C))
\end{equation}
and
\begin{equation}
H^0(X,\Ox_X(K_X+m(M+3H)+H)) \rightarrow H^0(X,\Ox_X(m(K_X+M+6H)-2H))
\end{equation}
are injective. Also, by Kawamata-Viehweg vanishing, we have 
\begin{equation}
H^1(X,\Ox_X(K_X+m(M+3H)))=0.
\end{equation}
Therefore, the map
\begin{equation}
H^0(X,\Ox_X(K_X+m(M+3H)+H)) \rightarrow H^0(C,\Ox_C(K_C+m(M_C+3H_C)))
\end{equation}
is surjective. Thus, by commutativity of the above diagram, we get
\begin{equation}
\begin{split}
h^0(\Ox_X(m(K_X+M+6H)-2H))-h^0(\Ox_X&(m(K_X+M+6H)-3H))=\\
\dim \mathrm{Im} (H^0(\Ox_X(m(K_X+M+&6H)-2H)) \rightarrow\\ \rightarrow H^0(\Ox_C(m(K_C+M_C+&5H_C)-2H_C))) \geq \\
h^0(\Ox_C(K_C+m(M_C+&3H_C))).
\end{split}
\end{equation}
Now, we notice that
\begin{equation}
m(K_X+M+6H)-3H=(m-1)(K_X+M+6H)+K_X+M+3H.
\end{equation}
Since $|K_X+M+3H|\neq 0$  and $H$ is effective, we get
\begin{equation} \label{inequality1}
\begin{split}
h^0(\Ox_X(m(K_X+M+6H)))&-h^0(\Ox_X((m-1)(K_X+M+6H)))\geq\\
h^0(\Ox_C(K_C&+m(M_C+3H_C))).
\end{split}
\end{equation}

Now, define $P(m) \coloneqq h^0(\Ox_X(m(K_X+M+6H)))$, and analogously set $Q(m) \coloneqq h^0(\Ox_C(K_C+m(M_C+3H_C)))$. Thus, inequality (\ref{inequality1}) can be rephrased as
\begin{equation}
P(m)-P(m-1)\geq Q(m).
\end{equation}
The difference $P(m)-P(m-1)$ is a linear polynomial in $m$ up to a periodic bounded perturbation \cite[p. 171]{LAZ1}, while $Q(m)$ is honestly linear in $m$. The two (almost) polynomials have leading coefficients $\frac{1}{2}\vol(X,K_X+M+6H)$ and $(M+3H)\cdot H$, respectively. Thus, we get
\begin{equation}
\frac{1}{2} \vol(X, K_X+M+6H) \geq (M+3H)\cdot H.
\end{equation}
Now, since $K_X+M\leq K_X+B+M$, and $H \leq m_0(K_X+B+M)$, we have
\begin{equation}
\begin{split}
H \cdot (M+3H) &\leq \frac{1}{2} \vol (X,(6m_0+1)(K_X+B+M))\\ &= \frac{(6m_0+1)^2}{2} \vol(X, K_X+B+M),
\end{split}
\end{equation}
and this concludes the proof. \hfill $\square$
\end{dimo1*}

The next statement takes care of the volume of the nef divisor $M$.

\begin{prop} \label{surfaceBoundM2}
Let $(X,B)$ be a log canonical pair, where $X$ is a smooth surface and $B \in \Lambda$, where $\Lambda \subset [0,1]$ is a DCC set. Also, let $M$ be a nef Cartier divisor. Assume that $K_X+B+M$ is big. Then, there is a constant $C$ such that
\begin{equation*}
M^2 \leq C \, \vol (X,K_X+B+M).
\end{equation*}
\end{prop}
\begin{dimo1*}
By \cite[Theorem 8.1]{BZ}, there is $e \in (0,1)$, depending just on $\Lambda$, such that $K_X+eB+eM$ is big. Write $K_X+eB+eM \sim_\R A +E$, where $A$ is ample and $E$ is effective. Then, we have
\begin{equation}
\begin{split}
\vol(X,K_X+B+M)&=\vol(X,K_X+eB+eM+(1-e)B+(1-e)M)\\
& \geq \vol(X,A+(1-e)M)\\
& \geq (1-e)^2 M^2.
\end{split}
\end{equation}
Thus, we can pick $C \coloneqq (1-e)^{-2}$. \hfill $\square$
\end{dimo1*}

\begin{rem}
Proposition \ref{surfaceBoundM2} goes through in the case $rM$ is Cartier for some positive integer $r$. Indeed, \cite[Theorem 8.1]{BZ} applies in case the coefficients of $M$ are in $\Lambda$; therefore, up to enlarging $\Lambda$ by adding $\N \frac{1}{r}$ to it, we can again apply \cite[Theorem 8.1]{BZ}.
\end{rem}

We conclude this section bounding the intersection product between the support of an element in a linear series as in Theorem \ref{BirkarZhangBirational} and the free part of the same linear series.

\begin{prop} \label{boundHG}
Let $(X,B)$ be a log canonical pair, where $X$ is a smooth surface and $B \in \Lambda$, where $\Lambda \subset [0,1]$ is a DCC set. Assume that $M$ is a nef $\R$-Cartier $\R$-divisor, and $|\lfloor m_0(K_X+B+M) \rfloor|=|H|+F$ defines a birational morphism. Then, we have the inequality
\begin{equation*}
G_\red \cdot H \leq \frac{2}{5}(14m_0+2)^2 \vol(X,K_X+B+M),
\end{equation*}
where $0 \leq G \sim \lfloor m_0(K_X+B+M) \rfloor$.
\end{prop}
\begin{dimo1*}
By \cite[Lemma 2.7.5]{HMX}, we have
\begin{equation}
\begin{split}
10 (G_\red \cdot H) & \leq 4 \, \vol(X,K_X+10H+G_\red)\\
& \leq 4 \, \vol(X,K_X+11m_0(K_X+B+M)).
\end{split}
\end{equation}
Now, as discussed in the proof of Proposition \ref{SurfaceBoundHM}, $|K_X+2M+3H|\neq \emptyset$. This, together with $B \geq 0$, provides us with
\begin{equation}
\begin{split}
10 (G_\red \cdot H) &\leq 4 \, \vol(X,K_X+11m_0(K_X+B+M)) \\ &\leq 4 \, \vol(X,(14m_0+2)(K_X+B+M)),
\end{split}
\end{equation}
and this concludes the proof. \hfill $\square$
\end{dimo1*}

\section{DCC for Volumes of Surfaces}

In this section, we prove the first main result, i.e. Theorem \ref{MainResult}. The proof strategy is as follows. For simplicity, we show the proof in case $r=1$, i.e. the nef part $M$ is Cartier. The general case reduces easily to this one. Then, we proceed by contradiction, and assume there is a sequence $\lbrace (X_i,B_i+M_i) \rbrace_{i \geq 1}$ as in the statement of Theorem \ref{MainResult} such that $\lbrace \vol(X_i,K_{X_i}+B_i+M_i) \rbrace_{i \geq 1}$ forms a strictly decreasing sequence. The idea is to put these log pairs in a family and use Theorem \ref{DefInvariance} and Theorem \ref{DCCcase1} to derive a contradiction.

Since the sequence of volumes is decreasing, it is bounded. By Theorem \ref{BirkarZhangBirational}, we may assume there is a positive integer $m_0$ such that $\lfloor m_0(K_{X_i}+B_i+M_i) \rfloor$ defines a birational morphism for all $i$. In view of Lemma \ref{smoothModel} and Proposition \ref{boundHB}, the log pairs $(X_i,B_i)$ are log birationally bounded. Thus, we obtain a bounding family $(\mathcal{X},\mathcal{B}) \rightarrow T$, which we may assume to be log smooth over the base. 

The main difficulty of the proof is to produce some divisor $\mathcal{M}$ on $\mathcal{X}$ in order to bound the $M_i$'s and be able to apply deformation invariance of volumes. In general, we have no control on the way a nef divisor is represented as linear combination of prime divisors. In particular, a nef divisor with given volume may have arbitrarily large coefficients, both in the positive and negative directions. Our strategy is to represent the $M_i$'s as a difference of divisors we have control upon. Since the linear series of $\lfloor m_0(K_{X_i}+B_i+M_i) \rfloor$ defines a birational morphism, we can pick an effective divisor $G_i \sim \lfloor m_0(K_{X_i}+B_i+M_i) \rfloor$. Thanks to Proposition \ref{boundHG}, we can bound the $G_i$'s as well. Therefore, we can write $m_0 M_i \sim G_i -m_0K_{X_i}-\lfloor m_0 B_i \rfloor$.

The last problem is to bound the components of $G_i$ that are exceptional for the morphism defined by $\lfloor m_0(K_{X_i}+B_i+M_i) \rfloor$. We can control them by relying on the extra information the surface assumption guarantees us. Once all of this is set up, and we have a divisor $\mathcal{M}$ bounding the $M_i$'s, we can proceed in the fashion of Theorem \ref{teoA}.

\begin{dimo1*}[{Theorem \ref{MainResult}, $r=1$ case}]
We will proceed by contradiction. If the statement does not hold, there is a sequence of generalized polarized pairs $(X_i,B_i+M_i) \in \mathfrak{V}$ such that the volumes $\vol(X_i,K_{X_i}+B_i+M_i)$ form a strictly decreasing sequence. By Theorem \ref{BirkarZhangBirational}, there is a positive integer $m_0$ such that $|m_0(K_{X_i}+B_i+M_i)|$ defines a birational map. By Lemma \ref{smoothModel}, we may assume the $X_i$'s are smooth, and that these maps are actual morphisms.

\underline{Step 1:} In this step we show that the sequence $\lbrace (X_i,B_i+M_i) \rbrace_{i \geq 1}$ is log birationally bounded.

Call $\phi_i:X_i \rightarrow Z_i$ the birational morphism defined by $|\lfloor m_0(K_{X_i}+B_i+M_i) \rfloor|$. By our assumptions, $\lfloor m_0(K_{X_i}+B_i+M_i)\rfloor \sim G_i = H_i + F_i \geq 0$, where $|H_i|$ defines $\phi_i$ and is a general element of such free linear series, and $F_i$ is the fixed part of $|\lfloor m_0(K_{X_i}+B_i+M_i)\rfloor|$. By Proposition \ref{boundHB} and Proposition \ref{boundHG}, the intersections $B_{i,\red} \cdot H_i$ and $G_{i,\red} \cdot H_i$ are bounded. Also, since the sequence $\lbrace \vol(X_i,K_{X_i}+B_i+M_i)\rbrace_{i \geq 1}$ is strictly decreasing, the self-intersections $H_i^2$ are bounded. Therefore, as explained in Remark \ref{BoundLogBir}, the pairs $(X_i,B_i+G_i)$ are log birationally bounded. Furthermore, by Proposition \ref{surfaceBoundM2} and Proposition \ref{SurfaceBoundHM}, the products $M_i^2$ and $M_i \cdot H_i$ are bounded as well.

By standard arguments, we can resolve the family that bounds the $Z_i$'s together with the $\phi_i(B_{i,\mathrm{red}})$'s and $\phi_i(G_{i,\mathrm{red}})$'s, and obtain a log smooth family $(\mathcal{X},\mathcal{B} +\mathcal{G}+\mathcal{E}) \rightarrow T$, where $T$ is smooth of finite type over $\C$, yet possibly reducible. Here, $\mathcal{B}$ and $\mathcal{G}$ are the strict transforms of the divisors that bound the birational transforms of the supports of the $B_i$'s and $G_i$'s, respectively, on the original singular family, while $\mathcal{E}$ is the exceptional divisor we introduce in resolving the family bounding the $Z_i$'s. We write $\mathcal{E}=\sum_{j=1}^l\mathcal{E}^j$ as sum of irreducible components. In particular, $\mathcal{B}+\mathcal{E}$ bounds the birational transforms of the $B_i$'s, while $\mathcal{G}+\mathcal{E}$ bounds the birational transforms of the $G_i$'s.

Since $T$ has finitely many components, up to passing to a subsequence, we may assume that all of the representatives of the $X_i$'s live over one fixed component of $T$. Therefore, up to replacing $T$, we may assume it is irreducible as well. By Lemma \ref{smoothModel}, we can assume that all the $X_i \dashrightarrow \mathcal{X}_i$ are morphisms, and that the divisors $B_i+F_i$ all have simple normal crossing support. Here $\mathcal{X}_i$ denotes the fiber of $\mathcal{X}\rightarrow T$ corresponding to $X_i$.

\underline{Step 2:} In this step, we show how to build a divisor that bounds the $M_i$'s.

On each $X_i$, we can write $m_0 M_i \sim G_i-m_0 K_{X_i}- \lfloor m_0 B_i \rfloor$. Analogously, since linear equivalence is preserved by pushforward, we have $m_0 M_i' \sim G_i'-m_0 K_{\mathcal{X}_i}-\lfloor m_0B_i'\rfloor$, where the prime notation denotes the pushforward to $\mathcal{X}_i$ of the corresponding divisors. To use Theorem \ref{DefInvariance}, we would need to have a divisor $\mathcal{M}$ on $\mathcal{X}$ that restricts (at least numerically) to the $M_i'$'s on all but finitely many $\mathcal{X}_i$'s.

Since we have a smooth family, $K_{\mathcal{X}/T}$ restricts fiberwise to the canonical sheaf of the fiber. Hence, $-m_0 K_{\mathcal{X}/T}$ restricts to $-m_0 K_{\mathcal{X}_i}$ for all $i$. Furthermore, by construction, $G_i'-\lfloor m_0B_i'\rfloor$ is supported on $\mathcal{B}_i+\mathcal{G}_i+\mathcal{E}_i$, where $\mathcal{B}_i+\mathcal{G}_i+\mathcal{E}_i \coloneqq (\mathcal{B}+\mathcal{G}+\mathcal{E})|_{\mathcal{X}_i}$. Therefore, we are left with bounding the coefficient of each irreducible component of $G_i'-\lfloor m_0B_i'\rfloor$. If we manage to do so, by the pigeonhole principle, there is a choice of coefficients for the components of $\mathcal{B}+\mathcal{G}+\mathcal{E}$ so that the corresponding divisor restricts to $G_i'-\lfloor m_0B_i'\rfloor$ on $\mathcal{X}_i$ for infinitely many $i$.

Now, since the $B_i$'s are effective with coefficients in $[0,1]$, the irreducible components of $\lfloor m_0 B_i'\rfloor$ have coefficients in $\lbrace 0 ,\ldots , m_0 \rbrace$. 
Therefore, we are left with showing that the coefficients of the $G_i'$'s live in a finite set as well. We notice that $G_i=H_i+F_i$, and $H_i$ is smooth and irreducible by Bertini's Theorem \cite[Corollary III.10.9]{Har77}. Thus, we just have to consider the coefficients of $F_i'$, which is the pushforward of $F_i$ onto $\mathcal{X}_i$. By Lemma \ref{SurfaceBoundHM}, $M_i \cdot H_i$ is bounded; hence, the coefficients of the components of $F_i$ that are not $\phi_i$-exceptional are bounded from above, since $H_i$ is the pullback of a very ample class on $Z_i$. Also, notice that $F_i$ is effective. Therefore, we are left with finding an upper bound for the coefficients of the $(\mathcal{X}_i \rightarrow Z_i)$-exceptional components of $F_i'$.

\underline{Step 3:} In this step, we bound the coefficients of the $(\mathcal{X}_i \rightarrow Z_i)$-exceptional components of the $F_i'$'s.

Denote by $Y_i$ the normalization of $Z_i$. Since $\mathcal{X}_i$ is smooth, $\mathcal{X}_i\rightarrow Z_i$ factors through $Y_i$, and this gives the Stein factorization of the morphism \cite[Corollary III.11.5]{Har77}. Call the induced map $f_i: \mathcal{X}_i \rightarrow Y_i$. By Lemma \ref{push nef},  $G_i'-m_0 K_{\mathcal{X}_i}-\lfloor m_0B_i'\rfloor$ is nef. By \cite[Lemma 3.40]{KM}, the intersection matrix of the $f_i$-exceptional curves is negative definite. Thus, for any Weil divisor on $Y_i$ we can define a numerical pullback on $\mathcal{X}_i$. Given a Weil divisor $D$, its numerical pullback is $f_{i,*}^{-1}D+\sum_{j=1}^l a_j \mathcal{E}_i^j$, where the coefficients $a_j$ are the unique solution to the system of equations
\begin{equation}
\left(f_{i,*}^{-1}D+\sum_{j=1}^l a_j \mathcal{E}_i^j\right)\cdot \mathcal{E}_i^k = 0, \quad 1 \leq k \leq l.
\end{equation}
The numerical pullback is defined so that the divisor obtained is $f_i$-trivial. By abusing notation, we will denote such an operator by $f_i^*$.

Now, assume $C$ is a nef divisor on $\mathcal{X}_i$. Then, $C-f_i^*f_{i,*}C$ is $f_i$-nef. Therefore, by the Negativity Lemma \cite[Lemma 3.39]{KM}, $C-f_i^*f_{i,*}C \leq 0$. Therefore, we have
\begin{equation}
G_i'-m_0 K_{\mathcal{X}_i}-\lfloor m_0B_i'\rfloor \leq f_i^*f_{i,*}(G_i'-m_0 K_{\mathcal{X}_i}-\lfloor m_0B_i'\rfloor).
\end{equation}

Now, by the above discussion, the coefficients of the components of $f_{i,*}^{-1}f_{i,*}(G_i'-m_0 K_{\mathcal{X}_i}-\lfloor m_0B_i'\rfloor)$ are bounded, since they correspond to the non-exceptional components of $G_i'$, $\lfloor m_0 B_i' \rfloor$ and $K_{\mathcal{X}_i}$. So, up to passing to a subsequence, we may assume they are constant as $i$ varies. Also, since all the components of such divisors are obtained restricting $K_{\mathcal{X}/T}$ and $\mathcal{B}+\mathcal{G}+\mathcal{E}$ to $\mathcal{X}_i$, we have a divisor $\mathcal{D}$ that on each $\mathcal{X}_i$ induces $f_{i,*}^{-1}f_{i,*}(G_i'-m_0 K_{\mathcal{X}_i}-\lfloor m_0B_i'\rfloor)$. Furthermore, the exceptional divisors for the maps $f_i$'s are supported on the divisors $\mathcal{E}_i$.

Now, fix $i_0$, and let $\mathcal{C}=\sum_{j=1}^l b_j \mathcal{E}^j$ be the divisor that restricts to the $f_{i_0}$-exceptional part of $f_{i_0}^* f_{{i_0},*}(G_{i_0}'-m_0 K_{\mathcal{X}_{i_0}}-\lfloor m_0B_{i_0}'\rfloor)$ on $\mathcal{X}_{i_0}$. Since $\mathcal{X}\rightarrow T$ is smooth, any line bundle on $\mathcal{X}$ is flat over $T$. Therefore, for any invertible sheaf $\mathcal{L}$ on $\mathcal{X}$, the Euler characteristic of $\mathcal{L}_t$ is independent of $t \in T$ \cite[Theorem III.9.9]{Har77}. Furthermore, the intersection product on smooth surfaces is defined as the sum of the Euler characteristics of certain line bundles \cite[Theorem I.4]{BEAU}. Hence, the intersection products
\begin{equation}
(\mathcal{D}+\mathcal{C})_t \cdot \mathcal{E}_t^j
\end{equation}
are independent of $t \in T$. Therefore, for each $i$, we have
\begin{equation} \label{inequalityZariskiPullback}
f_i^* f_{i,*}(G_i'-m_0 K_{\mathcal{X}_i}-\lfloor m_0B_{i}'\rfloor)=f_{i,*}^{-1} f_{i,*}(G_i'-m_0 K_{\mathcal{X}_i}-\lfloor m_0B_{i}'\rfloor) + \sum_{j=1}^l b_j \mathcal{E}_j^i.
\end{equation}

This argument shows that the coefficients of the $f_i$-exceptional components of the $f_i^* f_{i,*}(G_i'-m_0 K_{\mathcal{X}_i}-\lfloor m_0B_i'\rfloor)$'s are constant. By inequality (\ref{inequalityZariskiPullback}), the $f_i$-exceptional components of the $(G_i'-m_0 K_{\mathcal{X}_i}-\lfloor m_0B_i'\rfloor)$'s are bounded from above. Since the coefficients of the $f_i$-exceptional components of the $\lfloor m_0 B_i' \rfloor$'s and the $K_{\mathcal{X}_i}$'s are bounded, we conclude that the coefficients of the $f_i$-exceptional components of the $F_i'$'s are bounded from above as well.

\underline{Step 4:} In this step, we show that there is a divisor $m_0\mathcal{M}$ that, up to linear equivalence, restricts to $m_0 M_i'$ on infinitely many $\mathcal{X}_i$'s.

Since the coefficients of the $\lfloor m_0 B_i'\rfloor$'s and the $G_i'$'s are bounded, and their components are bounded by $\mathcal{B}+\mathcal{G}+\mathcal{E}$, up to passing to a subsequence, we may assume that there is a divisor $\mathcal{N}$ on $\mathcal{X}$ that cuts $G_i'-\lfloor m_0B_i'\rfloor$ on each $\mathcal{X}_i$. Now, we define $m_0 \mathcal{M} \coloneqq \mathcal{N}-m_0 K_{\mathcal{X}/T}$. The divisor $m_0\mathcal{M}$ is so that $m_0\mathcal{M}_i \sim m_0 M_i'$ for every $i$. 

Now, since there is a $\Q$-divisor $\mathcal{M}$ such that $\mathcal{M}_i \sim_\Q M_i'$ for all $i$, again, by flatness, we argue that $M_i'^2$ is independent of $i$. Since by Proposition \ref{surfaceBoundM2} the volumes $M_i^2$ are bounded, after passing to a subsequence, we may assume that $M_i^2$ is constant as well. So, the number $d \coloneqq M_i'^2-M_i^2$ is independent of $i$.

\underline{Step 5:} In this step, following the ideas in Theorem \ref{structureSurfaceNef}, we replace the family $\mathcal{X} \rightarrow T$ with a new family $\mathcal{X}' \rightarrow T'$ that parametrizes the blow-ups at $d$ points of the surfaces parametrized by $\mathcal{X}$. In this way, we can arrange for the $M_i$'s to descend to the new birational models.

We can produce such a new family inductively. To start, consider $\mathcal{X}\times_T \mathcal{X}$ with projections $p$ and $q$, and denote the diagonal over $T$ by $\Delta_T$. Then, we pull back along $q$ all the irreducible components of the relevant divisors on $\mathcal{X}$ (e.g. $\mathcal{B}$, $\mathcal{G}$ and $\mathcal{E}$). Notice that on $\mathcal{X}\times_T \mathcal{X}$ these meet transversally with $\Delta_T$. Then, we blow up $\Delta_T$. The family
\begin{equation}
\mathcal{X}^{(1)} \coloneqq Bl_{\Delta_T} \mathcal{X}\times_T \mathcal{X} \xrightarrow{p \circ bl_{\Delta_T}} T^{(1)} \coloneqq \mathcal{X}
\end{equation}
parametrizes the blow-ups at one point of the fibers of $\mathcal{X} \rightarrow T$. On this new family we have divisors $\mathcal{B}^{(1)}$, $\mathcal{G}^{(1)}$ and $\mathcal{E}^{(1)}$ that correspond to the strict transforms of $\mathcal{B}$, $\mathcal{G}$ and $\mathcal{E}$. Furthermore, we have the $bl_{\Delta_T}$-exceptional divisor $\mathcal{P}^{(1)}$. We iterate this process $d$ times, and we end up with divisors $\mathcal{B}^{(d)}$, $\mathcal{G}^{(d)}$, $\mathcal{E}^{(d)}$ and $\mathcal{P}^{(1)}, \ldots, \mathcal{P}^{(d)}$. Here, each $\mathcal{P}^{(k)}$ denotes the strict transform of the exceptional divisor obtained at the $k$-th step. The family $(\mathcal{X}^{(d)},\mathcal{B}^{(d)}+\mathcal{G}^{(d)}+\mathcal{E}^{(d)}+\sum_{k=1}^d \mathcal{P}^{(k)})\rightarrow T^{(d)}$ is in general not log smooth. Therefore, we further base change and resolve it (if necessary), and then decompose the base in finitely many locally closed subsets, so that each subfamily is log smooth over the base. Notice that in this process we may introduce further exceptional divisors, say $\mathcal{P}^{(d+1)},\ldots , \mathcal{P}^{(q)}$. By construction, we may assume that they are part of the log smooth family.

To summarize, we obtained a family
\begin{equation}
\left( \mathcal{X}',\mathcal{B}'+\mathcal{G}'+\mathcal{E}'+\sum_{k=1}^q \mathcal{P}^{(k)}\right)\rightarrow T',
\end{equation}
where $T'=\bigsqcup_{h=1}^p T'^{(h)}$ as a union of locally closed subsets, and the family is log smooth over each one of the $T'^{(h)}$. Furthermore, by Theorem \ref{structureSurfaceNef}, for each $i$ there is a fiber $\mathcal{X}'_i$ such that $f_i: X_i \rightarrow \mathcal{X}_i$ factors through $\mathcal{X}'_i$, and $M_i$ descends to $\mathcal{X}'_i$. Call $\overline{M}_i$ such divisor on $\mathcal{X}'_i$.

Now, define $\mathcal{M}'$ to be the pullback of $\mathcal{M}$ to $\mathcal{X}'$. As argued before, over each $T'^{(h)}$, the self-intersections $(\mathcal{M}'_t)^2$ and $(\mathcal{P}^{(k)}_t)^2$ are constant. Furthermore,  $(\mathcal{P}^{(k)}_t)^2<0$ for any $k$ and for any $T'^{(h)}$. Since $\mathcal{M}'_t \cdot \mathcal{P}^{(k)}_t=0$ for any $k$ and any $t \in T'$, there are finitely many $q$-tuples of non-negative integers $(c_1, \ldots , c_q)$ such that the divisor $\mathcal{M}'-\sum_{k=1}^qc_k \mathcal{P}^{(k)}$ restricts ($\Q$-linearly) to a $\overline{M}_i$ for some $i$. Therefore, by the pigeonhole principle, there is a choice of $(c_1,\ldots , c_q)$ such that $\overline{M}_i$ is the restriction of $(\mathcal{M}'-\sum_{k=1}^q c_k \mathcal{P}^{(k)})_{i,j(i)}$ for infinitely many $i$. Furthermore, we may assume that all the $\mathcal{X}'_{i}$ live over the same $T'^{(h)}$.

\underline{Step 6:} In this step, we use Theorem \ref{DefInvariance} and Theorem \ref{DCCcase1} to conclude the proof.

By the previous step, up to replacing $\mathcal{X} \rightarrow T$ with a locally closed subset of $\mathcal{X}' \rightarrow T'$, we may assume that $M_i'^2=M_i^2$, and that $M_i$ is obtained by pullback of $M_i'$ along $X_i \rightarrow \mathcal{X}_i$. Call such maps $g_i: X_i \rightarrow \mathcal{X}_i$. Now, by the proof of \cite[Corollary 3.1.3]{HMX} and Lemma \ref{volumes lemma}, we may assume that the $g_i$'s are induced by a finite sequence of blow-ups along strata of $(\mathcal{X}_i,\mathcal{B}_i)$. Since $(\mathcal{X},\mathcal{B})$ is log smooth, there is a sequence of blow-ups along strata of $\mathcal{B}$, say $\Psi^{(i)}: \overline{\mathcal{X}}^{(i)}\rightarrow \mathcal{X}$, such that $X_i \cong \overline{\mathcal{X}}_i^{(i)}$. Call $\overline{\mathcal{M}}^{(i)}$ the pullback of $\mathcal{M}$ to $\overline{\mathcal{X}}^{(i)}$. Also, let $\Gamma^{(i)}$ be the divisor supported on the strict transform of $\mathcal{B}$ and the $\Psi^{(i)}$-exceptional divisors such that $\Gamma_i^{(i)} = B_i$. Finally, fix some $i_0$. Then, we have
\begin{equation}
\begin{split}
\vol(X_i,K_{X_i}+B_i+M_i) &= \vol(\overline{\mathcal{X}}_{i}^{(i)},K_{\overline{\mathcal{X}}_{i}^{(i)}}+\Gamma_{i}^{(i)} + \overline{M}_{i}^{(i)}) \\ &= \vol(\overline{\mathcal{X}}_{i_0}^{(i)},K_{\overline{\mathcal{X}}_{i_0}^{(i)}}+\Gamma_{i_0}^{(i)} + \overline{M}_{i_0}^{(i)}),
\end{split}
\end{equation}
where the second equality follows from Theorem \ref{DefInvariance}. Thus, all the volumes in our sequence are realized by some model $(\overline{\mathcal{X}}_{i_0}^{(i)},\Gamma_{i_0}^{(i)} + \overline{M}_{i_0}^{(i)})$ living over $(\mathcal{X}_{i_0},B_{i_0}'+M_{i_0}')$. Therefore, by Theorem \ref{DCCcase1}, we get a contradiction. Hence, $\mathfrak{V}$ satisfies the DCC property. \hfill $\square$
\end{dimo1*}

\begin{rem}
Since all the facts needed for the proof go through in case the nef part is $\Q$-Cartier, the statement holds true in  case we allow the nef part to be $\Q$-Cartier with given Cartier index $r$. This recovers the full statement of Theorem \ref{MainResult}.
\end{rem}

\section{Towards Boundedness}

In this section, we discuss some results related to boundedness. First, we introduce some technical properties that hold in all dimensions. Then, we focus on the case of surfaces, and prove Theorem \ref{BoundednessLCModels}.

\subsection{Some General Facts about Boundedness}

The following is a generalization of \cite[Proposition 5.1]{HMX} in the setting of generalized polarized pairs.  It shows that, once we fix a base model $(Z,D)$ and a volume for the models living above $(Z,D)$, there exists a higher model $f:Z'\rightarrow Z$ where divisors naturally descend.

\begin{prop} \label{VolumeWgeneralCase}
Fix a positive number $w$, a non-negative integer $d$, and a DCC set $\Lambda \subset [0,1]$. Let $(Z,D)$ be a projective log smooth $d$ dimensional pair where $D$ is reduced, and let $M$ be a nef $\R$-Cartier $\R$-divisor on $Z$. Then, there exists $f:Z' \rightarrow Z$, a finite sequence of blow-ups along strata of the b-divisor $\mathbf{M}_D$, such that if
\begin{itemize}
\item $(X,B)$ is a projective log smooth $d$-dimensional pair;
\item $g : X \rightarrow Z$ is a finite sequence of blow-ups along strata of $\mathbf{M}_D$;
\item $\coeff (B) \subset \Lambda$;
\item $g_*B \leq D$;
\item $\vol (X,K_X+B+g^*M)=w$;
\end{itemize}
then $\vol(Z',K_{Z'}+\mathbf{M}_{B,Z'}+f^*M)=w$.
\end{prop}
\begin{dimo1*}
We can assume that $1 \in \Lambda$. Let $\mathfrak{P}$ be the set of pairs $(X,B)$ over $(Z,D)$ that satisfy the first four of the five hypotheses of the statement. Notice that $\mathfrak{P}$ is a subset of the set $\mathfrak{D}$ in Theorem \ref{DCCcase1}. Then, define
\begin{equation}
\mathfrak{W}  \coloneqq  \lbrace \vol(X,K_X+B+g^*M)|(X,B)\in \mathfrak{P} \rbrace.
\end{equation}
By Theorem \ref{DCCcase1}, $\mathfrak{W}$ satisfies the DCC. Therefore, there is a constant $\delta > 0$ such that, if $w \leq \vol(X,K_X+B+g^*M) \leq w + \delta$, then $\vol(X,K_X+B+g^*M)=w$. Also, by \cite[Theorem 8.1]{BZ}, there exists an integer $r$ such that, if $(X,B) \in \mathfrak{P}$ and $K_X+B+g^*M$ is big, then $K_X+\frac{r-1}{r}(B+g^*M)$ is big as well. Now, fix $\epsilon > 0$ such that $(1-\epsilon)^d> \frac{w}{w+\delta}$, and define $a  \coloneqq  1-\frac{\epsilon}{r}$.

Then, we have the following chain of inequalities
\begin{equation} \label{equation with fractions}
\begin{split}
\vol(X,K_X+a(B+g^*M)) &\geq \vol (X,(1-\epsilon)(K_X+B+g^*M))\\
&=(1-\epsilon)^d \vol(X,K_X+B+g^*M)\\
&> \frac{w}{w+\delta}\vol(X,K_X+B+g^*M),
\end{split}
\end{equation}
where the first inequality comes from the identity
\begin{equation}
K_X+a(B+g^*M)= (1-\epsilon)(K_X+B+g^*M)+\epsilon \left( K_X+\frac{r-1}{r}(B+g^*M)  \right) 
\end{equation}
and $K_X+\frac{r-1}{r}(B+g^*M)$ being big.

Now, since $(Z,aD)$ is Kawamata log terminal, we can obtain a terminalization $f:Z'\rightarrow Z$ by blowing up strata of $\mathbf{M}_D$. We can write
\begin{equation}
K_{Z'}+\Psi = f^*(K_Z+aD)+E,
\end{equation}
where $\Psi$ and $E$ are effective, $\Psi \wedge E =0$, and $(Z',\Psi)$ is terminal.

Let $\mathfrak{F}$ denote the set of pairs $(X,B)$ satisfying all the assumptions in the statement, and such that the rational map $\phi: X \dashrightarrow Z'$ is a morphism. Fix $(X,B) \in \mathfrak{F}$, and define $B_{Z'} \coloneqq \phi_*B$. Then, by construction, we have $f_*(aB_{Z'})\leq aD$. Therefore, if we write
\begin{equation}
K_{Z'}+\Phi= f^*(K_Z +f_* (aB_{Z'}))+F,
\end{equation}
where $\Phi$ and $F$ are effective with $\Phi \wedge F =0$, then $(Z',\Phi)$ is terminal. Hence, it follows that
\begin{equation} \label{equationTerminalization1}
\begin{split}
\vol(Z',K_{Z'}+aB_{Z'}+af^*M) &= \vol(Z',K_{Z'}+aB_{Z'}\wedge \Phi +af^*M)\\
&= \vol(X,K_{X}+\phi_*^{-1}(aB_{Z'}\wedge \Phi) +ag^*M)\\
& \leq \vol(X,K_{X}+B+g^*M),
\end{split}
\end{equation}
where the first equality follows from part (3) of Lemma \ref{volumes lemma}, the second one from part (2) of Lemma \ref{volumes lemma} and $(Z',aB_{Z'}\wedge \Phi)$ being terminal, and the last inequality from $\phi_*^{-1}(aB_{Z'}\wedge \Phi) \leq B$ and $M$ being nef.

Then, we get the following chain of inequalities
\begin{equation} \label{equationTerminalization2}
w \leq \vol(Z',K_{Z'}+B_{Z'}+f^*M) \leq \frac{w+ \delta}{w} \vol(Z',K_{Z'}+a(B_{Z'}+f^*M)) \leq w+\delta,
\end{equation}
where the first one follows from part (1) of Lemma \ref{volumes lemma}, the second one from inequality (\ref{equation with fractions}), and the third one from inequality (\ref{equationTerminalization1}). Therefore, by definition of $\delta$, we have $\vol(Z',K_{Z'}+B_{Z'}+f^*M)=w$.

To conclude the proof, it is enough to notice that, if $(X,B)$ satisfies the assumptions in the statement, then after blowing up along finitely many strata of $\mathbf{M}_D$ and replacing $B$ by its strict transform plus the exceptional divisors, we may assume that $(X,B) \in \mathfrak{F}$. \hfill $\square$
\end{dimo1*}

The following is a generalization of \cite[Lemma 5.3]{HMX}, and allows us to relate the generalized log canonical models of pairs with comparable boundary and same volume.

\begin{lemma} \label{LemmaSameLCModel}
Let $(X,B+M)$ be a generalized polarized pair such that $M$ is $\R$-Cartier, and $X=X'$, $M=M'$, and $(X,B)$ is log canonical. Assume that $(X,B+M)$ has a generalized log canonical model $f: X \dashrightarrow W$. If $B' \geq B$ is such that $(X,B')$ is log canonical, and $\vol(X,K_X+B+M)=\vol(X,K_X+B'+M)$, then $f$ is also the generalized log canonical model of $(X,B'+M)$.
\end{lemma}

\begin{dimo1*}
Up to going to a higher model\footnote{In doing so, as boundary we take the strict transforms of of $B$ (respectively $B'$) plus the reduced exceptional divisors, while as nef part we pull back $M$.}, we may assume that $f: X \rightarrow W$ is a morphism. Define $A \coloneqq  f_* (K_X+B+M)$. Then, $A$ is ample, and, as $W$ is a generalized log canonical model, $K_X+B+M-f^*A$ is effective and $f$-exceptional. Then, for any $t \in [0,1]$ we have
\begin{equation}
\begin{split}
\vol(X,K_X+B+M) &= \vol(X,K_X+B+t(B'-B)+M)\\
& \geq \vol(X,f^*A+t(B'-B)) \\
& \geq \vol(X,f^*A)\\
&= \vol(X,K_X+B+M).
\end{split}
\end{equation}

Therefore, by the proof of \cite[Lemma 5.3]{HMX}, $E \coloneqq B'-B$ is $f$-exceptional. Hence, we have
\begin{equation}
H^0(X,m(K_X+B'+M))=H^0(X,m(K_X+B+M))=H^0(W,mA)
\end{equation}
for all $m \geq 0$. In particular, this shows that $f: X \rightarrow W$ is the generalized log canonical model of $(X,B'+M)$. \hfill $\square$
\end{dimo1*}

\subsection{Boundedness for Surfaces}

Now, we are ready to prove Theorem \ref{BoundednessLCModels}. The proof goes as follows. By Theorem \ref{MainResult}, we can reduce to the case when our surfaces are irreducible. Arguing by contradiction, there is counterexample sequence $\lbrace (X_i,B_i+M_i)\rbrace_{i \geq 1}$ for which there is no uniform $N$ such that $N(K_{X_i}+B_i+M_i)$ is very ample for all $i$. Then, following the ideas in the proof of Theorem \ref{MainResult}, we can show that the sequence $\lbrace (X_i,B_i+M_i)\rbrace_{i \geq 1}$ is log birationally bounded by a log smooth family $(\mathcal{X},\mathcal{B}+\mathcal{M}) \rightarrow T$.

Using Proposition \ref{VolumeWgeneralCase} and Lemma \ref{LemmaSameLCModel}, we can reduce to the case when the maps $\mathcal{X}_i \dashrightarrow X_i$ are morphisms, and $X_i$ is the generalized log canonical model of $(\mathcal{X}_i,\mathbf{M}_{B_i,\mathcal{X}_i}+\mathcal{M}_i)$. Analyzing the morphisms $q_i:\mathcal{X}_i \rightarrow X_i$, we are able to argue that there is a unique divisor $\Phi$ supported on $\mathcal{B}$ that restricts to $\mathbf{M}_{B_i,\mathcal{X}_i}$ for infinitely many $i$.

Therefore, up to passing to a subsequence, we are reduced to the case where the coefficients of the $B_i$'s are fixed, and the $q_i$-exceptional divisors deform in family as a fixed subset of the components of $\mathcal{B}$. In particular, this shows that the Cartier index of the $X_i$'s is the same. 
This guarantees the existence of a positive integer $l$ such that $l(K_{X_i}+B_i+M_i)$ is integral and Cartier for all $i$.  Since $K_{X_i}+B_i+M_i$ is ample and $M_i$ is nef, we can then apply an effective basepoint-free theorem for semi-log canonical surfaces by Fujino \cite{Fuj17}. More precisely, we show that $12l(K_{X_i}+B_i+M_i)$ is very ample, and  the required contradiction follows.

\begin{dimo1*}[{Theorem \ref{BoundednessLCModels}, $r=1$ case}]
First, we notice that $w=\sum w_j$, where $w_j \coloneqq (K_{X_j}+B_j+M_j)^2$. By Theorem \ref{MainResult}, the $w_j$'s belong to a DCC set. Therefore, there are just finitely many admissible values for them. Thus, we can reduce to the case when $X$ is irreducible. Then, by Remark \ref{BoundLogBir}, it is enough to find an integer $N>0$ such that for any generalized polarized log canonical surface pair $(X,B+M)$ with $M$ is nef and Cartier, $\coeff(B) \in \Lambda$ and $(K_X+B+M)^2=w$, then $N(K_X+B+M)$ is very ample\footnote{Notice that in this formulation there is a hidden part of the claim: Once we fix the volume $w$, the denominators allowed in $\coeff(B)$ are bounded. This is implicit in $N(K_X+B+M)$ being Cartier.}. Also, without loss of generality, we may assume $1 \in \Lambda$.

\underline{Step 1:} Proceeding by contradiction, in this step, we build a log birationally bounded sequence of counterexamples that satisfies certain natural properties.

If the claim is not true, then there exists a sequence of irreducible surfaces $\lbrace (X_i,B_i+M_i)\rbrace_{i \geq 1} \subset \mathfrak{F}_{glc}(2,w,\Lambda)$ such that $i!(K_{X_i}+B_i+M_i)$ is not very ample for all $i \geq 1$. By the proof of Theorem \ref{MainResult}, we can bound the sequence of counterexamples at least log birationally. We get a log smooth family $(\mathcal{X},\mathcal{B}+\mathcal{M}) \rightarrow T$ such that for every $i$ we have a commutative diagram of birational maps

\begin{center}
\begin{tikzpicture}
  \matrix (m) [matrix of math nodes,row sep=3em,column sep=4em,minimum width=2em]
  {
     (Y_i,\mathbf{M}_{B_i,Y_i} + h_i^* M_i) & &  \\
      & (\mathcal{X}_i,\mathbf{M}_{B_i,\mathcal{X}_i}+\mathcal{M}_i) & \\
     (X_i,B_i+M_i) & & Z_i \\};
  \path[-stealth]
    (m-1-1) edge node [left] {$h_i$} (m-3-1)
            edge node [above right] {$g_i$} (m-2-2)
    (m-3-1) edge[dashed] node [below] {$\phi_{i}$} (m-3-3)
    (m-2-2) edge node [above right] {$f_{i}$} (m-3-3);
\end{tikzpicture}
\end{center}
where $\phi_i$ is induced by $| \lfloor m_0(K_{X_i}+B_i+M_i) \rfloor |$ for a fixed $m_0$ as in Theorem \ref{BirkarZhangBirational}, and $f_i\circ g_i$ is a resolution of $\phi_i$, which we can assume to be factoring through $\mathcal{X}_i$ by Lemma \ref{smoothModel}. Furthermore, by the proof of Theorem \ref{MainResult}, $m_0 \mathcal{M}_i$ is nef Cartier, and $g_i^*\mathcal{M}_i \sim_\Q h_i^*M_i$.

Define $\Phi_i^{(i)}  \coloneqq  \mathbf{M}_{B_i,\mathcal{X}_i}$, and let $\Phi^{(i)}$ be the divisor supported on $\mathcal{B}$ that restricts to $\Phi_i^{(i)}$ on $\mathcal{X}_i$. Since $\Lambda$ is a DCC set, up to passing to a subsequence, we may assume $\Phi^{(i)} \leq \Phi^{(i+i)}$ for all $i \geq 1$. Furthermore, as in the proof of Theorem \ref{MainResult}, we may assume that $M_i^2$ is constant as $i$ varies (this choice is already implicit in $h_i^*M_i \sim_\Q g_i^*\mathcal{M}_i$).

Now, let $\mathcal{W}_1 \rightarrow \mathcal{X}_1$ be the sequence of blow-ups along strata of $\mathbf{M}_{\mathcal{B}_1}$ constructed in Proposition \ref{VolumeWgeneralCase}. Since $(\mathcal{X},\mathcal{B})$ is log smooth, we can assume that $\mathcal{W}_1$ appears as fiber of a sequence of blow-ups $p: \mathcal{W}\rightarrow \mathcal{X}$ along strata of $\mathbf{M}_{\mathcal{B}}$. Therefore, up to replacing $(\mathcal{X},\mathcal{B}+\mathcal{M})$ with $(\mathcal{W},\mathbf{M}_{\mathcal{B},\mathcal{W}}+p^*\mathcal{M})$ and the $Y_i$'s with a higher model, we may assume that $(\mathcal{X}_{1},\mathcal{B}_{1}+\mathcal{M}_{1})$ satisfies the special property in Proposition \ref{VolumeWgeneralCase}. In particular, given any higher model $(\tilde{\mathcal{X}}_1,\tilde{B}) \xrightarrow{\alpha} \mathcal{X}_1$ obtained by blow-ups along strata of $\mathbf{M}_{\mathcal{B}_1}$ and such that $\alpha_ * \tilde{B} \leq \mathcal{B}_1$, $\tilde{B} \in \Lambda$ and $\vol(\tilde{\mathcal{X}}_1,K_{\tilde{\mathcal{X}}_1}+\tilde{B}+\alpha^*\mathcal{M}_1)=w$, then $\vol(\mathcal{X}_1,K_{\mathcal{X}_1}+\mathbf{M}_{\tilde{B},\mathcal{X}_1}+\mathcal{M}_1)=w$.

\underline{Step 2:} In this step, we show that $\vol(\mathcal{X}_i,K_{\mathcal{X}_i}+\Phi^{(i)}_i+\mathcal{M}_i)=w$ for all $i$.

By construction, this is true for $i=1$. Fix some $i \geq 2$. Then, let $\mathcal{E}_i$ be the finite set of divisors $E$ on $Y_i$ such that $a_E(\mathcal{X}_i,\mathbf{M}_{B_i,\mathcal{X}_i})<0$. Then, there is a finite sequence of blow-ups $\pi^{(i)}: \mathcal{X}^{(i)} \rightarrow \mathcal{X}$ along strata of $\mathbf{M}_{\mathcal{B}}$ such that the divisors in $\mathcal{E}_i$ are not exceptional for $Y_i \dashrightarrow \mathcal{X}^{(i)}_i$. Let $\Psi^{(i)}$ be the unique divisor supported on $\mathbf{M}_{\mathcal{B},\mathcal{X}^{(i)}}$ that restricts to $\mathbf{M}_{B_i,\mathcal{X}_i^{(i)}}$. Then, by part (4) of Lemma \ref{volumes lemma}, we have
\begin{equation}
\vol(\mathcal{X}_i^{(i)},K_{\mathcal{X}^{(i)}_i}+\Psi^{(i)}_i+\pi_i^{(i),*}\mathcal{M}_i)= \vol(Y_i,K_{Y_i}+\mathbf{M}_{B_i,Y_i}+g_i^*\mathcal{M}_i)=w,
\end{equation}
where $\pi_i^{(i)}$ denotes the restriction of $\pi^{(i)}$ to $\mathcal{X}^{(i)}_i$. Also, notice that, by construction, $\mathbf{M}_{\Psi^{(i)}_1,\mathcal{X}_1}=\Phi^{(i)}_1$. Therefore, this leads to the chain of equalities
\begin{equation}
\begin{split}
w &= \vol(\mathcal{X}_i^{(i)},K_{\mathcal{X}^{(i)}_i}+\Psi^{(i)}_i+\pi_i^{(i),*}\mathcal{M}_i) \\
&= \vol(\mathcal{X}_1^{(i)},K_{\mathcal{X}^{(i)}_1}+\Psi^{(i)}_1+\pi_1^{(i),*}\mathcal{M}_1) \\
&= \vol(\mathcal{X}_1,K_{\mathcal{X}_1}+\Phi_1^{(i)}+\mathcal{M}_1) \\
&= \vol(\mathcal{X}_i,K_{\mathcal{X}_i}+\Phi_i^{(i)}+\mathcal{M}_i),
\end{split}
\end{equation}
where the second and fourth ones follow from Theorem \ref{DefInvariance}, while the third one comes from the end of Step 1.

Also, by Lemma \ref{LemmaSameLCModel}, this computation shows that $(Y_i,\mathbf{M}_{\mathbf{M}_{B_i,\mathcal{X}_i},Y_i}+g_i^*\mathcal{M}_i)$ has the same generalized log canonical model as $(Y_i,\mathbf{M}_{B_i,Y_i} + h_i^* M_i)$, namely $(X_i,B_i+M_i)$. Furthermore, the former has the same section ring as $(\mathcal{X}_i,\mathbf{M}_{B_i,\mathcal{X}_i}+\mathcal{M}_i)$. Therefore, $(X_i,B_i+M_i)$ is the generalized log canonical model of $(\mathcal{X}_i,\mathbf{M}_{B_i,\mathcal{X}_i}+\mathcal{M}_i)$. In particular, since we are dealing with surfaces, the map from $\mathcal{X}_i$ to $X_i$ is an actual morphism. Thus, the above diagram of maps can be completed as follows

\begin{center}
\begin{tikzpicture}
  \matrix (m) [matrix of math nodes,row sep=3em,column sep=4em,minimum width=2em]
  {
     (Y_i,\mathbf{M}_{B_i,Y_i} + h_i^* M_i) & &  \\
      & (\mathcal{X}_i,\mathbf{M}_{B_i,\mathcal{X}_i}+\mathcal{M}_i) & \\
     (X_i,B_i+M_i) & & Z_i \\};
  \path[-stealth]
    (m-1-1) edge node [left] {$h_i$} (m-3-1)
            edge node [above right] {$g_i$} (m-2-2)
    (m-3-1) edge[dashed] node [below] {$\phi_{i}$} (m-3-3)
    (m-2-2) edge node [above right] {$f_{i}$} (m-3-3)
            edge node [above left] {$q_{i}$} (m-3-1);
\end{tikzpicture}
\end{center}

\underline{Step 3:} In this step, we show that the coefficients of the $B_i$'s live in a finite set.

By construction, we have $0 \leq \Phi^{(i)} \leq \Phi^{(i+1)} \leq \mathcal{B}$ for all $i \geq 1$. Therefore, we can define $\Phi^{(\infty)} \coloneqq \lim \Phi^{(i)}$. The previous step shows that for any $i$ and $j$ we have $\vol(\mathcal{X}_i,K_{\mathcal{X}_i}+\Phi^{(j)}_i+\mathcal{M}_i)=w$. Therefore, by continuity of the volume function, for all $i$ we have $\vol(\mathcal{X}_i,K_{\mathcal{X}_i}+\Phi^{(\infty)}_i+\mathcal{M}_i)=w$.

By Lemma \ref{LemmaSameLCModel}, $X_i$ is the generalized log canonical model of $(\mathcal{X}_i,K_{\mathcal{X}_i}+\Phi^{(\infty)}_i+\mathcal{M}_i)$, and the difference $\Phi^{(\infty)}_i-\Phi^{(i)}_i$ is $q_i$-exceptional. Now, let $\mathcal{P}$ be an irreducible component of $\mathcal{B}$. If $\mathcal{P}_i$ is $q_i$-exceptional for all $i \geq 1$, then $\mult_\mathcal{P} \Phi^{(i)}=1$ for all $i$. Thus, we have $\mult_\mathcal{P}\Phi^{(\infty)}=1$. On the other hand, if $\mathcal{P}_i$ is not $q_i$-exceptional for some $i$, the proof of Lemma \ref{LemmaSameLCModel} shows that $\mult_\mathcal{P}\Phi^{(\infty)}=\mult_\mathcal{P}\Phi^{(i)}$, and this is the coefficient of an irreducible component of $B_i$. Since $\Phi^{(\infty)}$ is supported on $\mathcal{B}$, it has only finitely many irreducible components. This shows that the coefficients of the $B_i$'s live in a finite set. Therefore, there exists $n>0$ such that $n(K_{X_i}+B_i)$ is integral Weil for all $i \geq 1$.

\underline{Step 4:} In this step, we show that the $X_i$'s all have the same Cartier index.

By definition of $\mathcal{B}$, up to passing to a subsequence, we may assume that any irreducible component $\mathcal{P}$ of $\mathcal{B}$ is either $q_i$-exceptional for all $i\geq 1$ or for no $i$. Furthermore, by the previous step, we may assume that $\Phi^{(i)}=\Phi^{(j)}$ for all $i,j \geq 1$.

In particular, the graphs of the $q_i$-exceptional curves and the corresponding weights are independent of $i$. Thus, by \cite[Remark 4.9]{KM}, we know that the singularities of the $X_i$'s are analytically isomorphic. Therefore, they have the same Cartier index.

\underline{Step 5:} In this step, we conclude the proof.

By the previous two steps, there exists a positive integer $l$ such that $l(K_{X_i}+B_i+M_i)$ is Cartier for all $i$. Since $K_{X_i}+B_i+M_i$ is ample and $M_i$ is nef, then the divisor $A_i  \coloneqq  4l(K_{X_i}+B_i+M_i)-(K_{X_i}+B_i)$ is ample. Furthermore, we have
\begin{equation}
A_i^2 \geq (3l(K_{X_i}+B_i+M_i))^2 \geq 9
\end{equation}
and
\begin{equation}
A_i \cdot C \geq 3l(K_{X_i}+B_i+M_i) \cdot C \geq 3
\end{equation}
for every curve $C$ in $X_i$. Therefore, by \cite[Theorem 5.1]{Fuj17}, the linear series $|4l(K_{X_i}+B_i+M_i)|$ is basepoint-free. Thus, by \cite[Lemma 7.1]{Fuj17}, $12l(K_{X_i}+B_i+M_i)$ is very ample. Hence, there is a uniform $N$ such that $N(K_{X_i}+B_i+M_i)$ is very ample for all $i$. This provides the required contradiction and concludes the proof. \hfill $\square$
\end{dimo1*}

\begin{rem}
As for the case of Theorem \ref{MainResult}, the proof of Theorem \ref{BoundednessLCModels} goes through with minor changes in case we fix the Cartier index of the nef part $M$.
\end{rem}

Now, we are ready to prove Theorem \ref{BigResult}. The statement follows from Theorem \ref{BoundLogBir} via an application of an effective basepoint-free theorem for semi-log canonical surfaces by Fujino \cite{Fuj17}.

\begin{dimo1*}[{Theorem \ref{BigResult}}]
Let $(X,B+M) \in \mathfrak{F}_{gslc}(2,w,\Lambda,r)$ and denote by $\psi: X^\nu \rightarrow X$ its normalization. Then, we can decompose $(X^\nu,\psi^*B+\psi^*M) \in \mathfrak{F}_{glc}(2,w,\Lambda,r)$ as a disjoint union of generalized log canonical pairs $(X_i,B_i+M_i)$. Furthermore, by Theorem \ref{BoundednessLCModels}, the generalized polarized pairs $(X_i,B_i+M_i)$ are bounded. In particular, there is a positive integer $N$ such that $N(K_{X_i}+B_i+M_i)$ is very ample.

Now, consider $A \coloneqq 5N(K_X+B+M)-(K_X+B)$. Since $K_X+B+M$ is ample and $M$ is nef, $A$ is ample. Also, let $Y_i$ be the irreducible component of $X$ whose normalization is $X_i$. Then, by construction, we have
\begin{equation}
A_i^2 = (\psi^*A_i)^2 \geq (4N(K_{X_i}+B_i+M_i))^2 \geq 16,
\end{equation}
where $A_i$ denotes the restriction of $A$ to $Y_i$. Furthermore, let $C \subset X$ be an irreducible curve. Then $C \subset Y_i$ for some $i$. Also, we can represent $A$ in a way so that $A \cdot C = A_i \cdot C$ is given by point counting. Furthermore, we may assume that these points lie in the smooth or nodal locus of $X$. Therefore, the normalization is at worst $2:1$ over the intersection points. Thus, we have
\begin{equation}
A \cdot C = A_i \cdot C \geq \frac{1}{2} (\psi^*A_i \cdot \psi^{-1}_* C) \geq 2N(K_{X_i}+B_i+M_i) \cdot \psi^{-1}_* C \geq 2.
\end{equation}
Hence, by \cite[Theorem 5.1]{Fuj17}, $|5N(K_X+B+M)|$ is basepoint-free. This, together with \cite[Lemma 7.1]{Fuj17}, implies that $15N(K_X+B+M)$ is very ample. Therefore there is a fixed multiple of $K_X+B+M$ that is very ample for all $(X,B+M) \in \mathfrak{F}_{gslc}(2,w,\Lambda,r)$, and the claim follows. \hfill $\square$
\end{dimo1*}

\begin{rem}
Since the boundedness in Theorem \ref{BigResult} is realized via a common multiple of $K_X+B+M$, and $(X,B)$ is log bounded, we conclude that $M$ is bounded up to $\Q$-linear equivalence as well.
\end{rem}

\section{Some Relevant Examples} \label{counterexamples}

In this section, we collect some examples of bad behaviors of generalized polarized pairs.

\subsection{Failure of Deformation Invariance of Plurigenera}

Let $E$ be an elliptic curve embedded as projectively normal subvariety of some projective space $\Pro^N$ with degree $n$, and let $X$ be the projective cone over it. Then, blowing up the vertex $v \in X$, we get a resolution $\pi : Y \rightarrow X$, where $Y$ is a $\Pro^1$-bundle over $E$, and $E_v \coloneqq \pi^{-1}(v)$ is a copy of $E$ with $E_v^2=-n$. It is well known that $X$ is log canonical but not Kawamata log terminal. In particular, we have
\begin{equation}
K_Y+E_v=\pi^*K_X.
\end{equation}
Denote by $L$ the class of a line through $v \in X$. Then, we have $K_X \equiv -nL$, and $nL \sim H$, where $H$ is very ample. Define $M \coloneqq 2H$. Then, $K_X+M$ is very ample. Therefore, $K_Y+E_v+\pi^*M=\pi^*(K_X+M)$ is nef and big. Let $g: Y \rightarrow E$ denote the $\Pro^1$-fibration. In particular, since $E_v$ is a section of $g$, for any line bundle $\mathcal{L}$ on $E$, the line bundle $g^*\mathcal{L}$ restricts to $\mathcal{L}$ on $E_v$.

Now, let $\mathcal{P}$ be the Poincar\'{e} line bundle on $E\du \times E$. Then, define the line bundle $\mathcal{Q} \coloneqq (\mathrm{id}_{E\du} \times g)^*\mathcal{P}$. Also, notice that the projection map $p: E\du \times Y \rightarrow E\du$ is smooth, and, in particular, flat. Call $q$ the other projection map $q: E\du \times Y \rightarrow Y$. On $E\du \times Y$, we have the line bundle $\mathcal{Q} \otimes q^*(K_Y+E+\pi^*M)$. By flatness of $p$, the Hilbert polynomial of the line bundles $(\mathcal{Q} \otimes q^*(K_Y+E+\pi^*M))|_{Y_P}$ is independent of $P \in E\du$ \cite[Theorem III.9.9]{Har77}, where $Y_P \coloneqq \lbrace P \rbrace \times Y$ and $\mathcal{Q}|_{Y_P}=g^* P$.

Multiplication by $E_v$ leads to the following short exact sequence
\begin{equation}
0 \rightarrow \Ox_{Y_P}(m(A+P')-E_{v}) \xrightarrow{\cdot E_v} \Ox_{Y_P}(m(A+P')) \rightarrow \Ox_{E_v}(mP) \rightarrow 0
\end{equation}
for any integer $m \geq 1$, where $A \coloneqq \pi^*(K_X+M)$ and $P' \coloneqq g^*P$. Since $A$ is the pullback of a very ample line bundle, when $P=0$ this leads to the short exact sequence of groups
\begin{equation} \label{SEScase0}
0 \rightarrow H^0(Y_0,\Ox_{Y_0}(mA-E_{v})) \rightarrow H^0(Y_0,\Ox_{Y_0}(mA)) \rightarrow \C \rightarrow 0.
\end{equation}

On the other hand, if $P$ is a non-torsion element of $E\du$, we have the vanishing $H^0(E_v,\Ox_{E_v}(mP))=0$ for all $m \geq 1$. This implies the following isomorphism
\begin{equation} \label{SEScaseP}
H^0({Y_P},\Ox_{Y_P}(m(A+P')-E_{v})) \cong H^0(Y_P,\Ox_{Y_P}(m(A+P')))
\end{equation}
for all $m \geq 1$.

In case $0 \neq P \in E\du$ is a torsion element, global sections behave as in equations (\ref{SEScase0}) or (\ref{SEScaseP}) depending on whether the order of $P$ divides $m$ or not.

By upper semi-continuity of the dimension of the space of global sections \cite[Theorem III.12.8]{Har77}, we get
\begin{equation}
h^0(Y_P,\Ox_{Y_P}(m(A+P'))) +1 \leq h^0(Y_0,\Ox_{Y_0}(mA)),
\end{equation}
where $P$ is not an $m$-torsion element of $E\du$. 

By construction, the family $(E\du \times Y,E\du \times E_v) \rightarrow E\du$ is log smooth. Also, the divisor $\mathcal{M} \coloneqq \mathcal{Q} \otimes (\pi \circ q)^* M$ restricts to a nef divisor on any fiber of the family. In particular, the divisor $K_{E\du \times Y/E\du}+E\du \times E_v+\mathcal{M}$ restricts to a nef and big divisor on every fiber. On the other hand, the example just constructed shows that deformation invariance of plurigenera does not hold for $K_{E\du \times Y/E\du}+E\du \times E_v+\mathcal{M}$. Therefore, deformation invariance of volumes, proved in Theorem \ref{DefInvariance}, is the best we can achieve in the setting of generalized polarized pairs.

\subsection{Failure of Finite Generation for the Pluricanonical Ring}

Throughout this section, we will refer to the notation just introduced in the previous example. We will show that the pluricanonical ring is not in general finitely generated for generalized polarized pairs.

Consider the nef and big divisor $\pi^*(K_X+M)+P'=K_Y+E+P'$ on $Y$, where $P'=g^*P$ and $P$ is a non-torsion element in $E\du$. Then, equation (\ref{SEScaseP}) shows that $\pi^*(K_X+M)+P'$ is not semiample. Therefore, by \cite[Theorem 2.3.15]{LAZ1}, the section ring $R(Y,\pi^*(K_X+M)+P')$ is not finitely generated.

\addcontentsline{toc}{section}{\bibname}
\printbibliography

\Addresses

\end{document}